\documentclass[11pt]{amsart}
\title{Weights on bimodules}
\author[\large P\lowercase{aramita} D\lowercase{as and} S\lowercase{hamindra} G\lowercase{hosh}]
{\bf \large P\lowercase{aramita} D\lowercase{as and} S\lowercase{hamindra} K\lowercase{umar} G\lowercase{hosh} }

\thanks{
Das was supported by K. U. Leuven grant BOF Research Grant OT/08/032 and Ghosh was supported by ERC Starting Grant VNALG-200749
}

\date{}

\address{Department of Mathematics, Katholieke Universiteit Leuven, Belgium}
\email{paramita.das@wis.kuleuven.be, shami.Ghosh@wis.kuleuven.be}

\usepackage{amsmath}
\usepackage{amsfonts}
\usepackage{latexsym}
\usepackage{amssymb}
\usepackage{mathrsfs}
\usepackage{amscd}
\usepackage{hyperref}
\usepackage{psfrag}
\usepackage{graphicx}
\usepackage{fullpage}
\usepackage[all]{xy}
\usepackage{verbatim}
\usepackage{enumerate}
\newtheorem{thm}{Theorem}[section]
\newtheorem{lem}[thm]{Lemma}

\newtheorem{prop}[thm]{Proposition}
\newtheorem{defn}[thm]{Definition}
\theoremstyle{remark}\newtheorem{rem}[thm]{Remark}
\theoremstyle{definition}\newtheorem{example}[thm]{Example}

\newenvironment{pf}{\noindent{\em Proof:}}{\hfill $\Box$\vspace*{1mm}}
\numberwithin{equation}{section}
\numberwithin{figure}{section}

\newcommand{\comments}[1]{}
\newcommand{\ra}{\rightarrow}
\newcommand{\lra}{\longrightarrow}

\newcommand{\mbb}{\mathbb}
\newcommand{\mcal}{\mathcal}

\newcommand{\N}{\mathbb N}
\newcommand{\Z}{\mathbb Z}
\newcommand{\C}{\mathbb{C}}

\newcommand{\mscr}{\mathscr}
\newcommand{\vlon}{\varepsilon}

\newcommand{\oline}{\overline}
\newcommand{\uline}{\underline}
\newcommand{\vphi}{\varphi}

\newcommand{\lsub}[2]{{\vphantom{#2}}_{#1}{#2}}


\begin{document}
\global\long\def\vlon{\varepsilon}
\global\long\def\ul#1{\underline{#1}}
\global\long\def\ol#1{\overline{#1}}
\global\long\def\os#1#2{\overset{#1}{#2}}
\global\long\def\us#1#2{\underset{#1}{#2}}
\global\long\def\ous#1#2#3{\overset{#1}{\underset{#3}{#2}}}
\global\long\def\lab{\langle}
\global\long\def\rab{\rangle}
\global\long\def\lrsuf#1#2#3{\vphantom{#2}_{#1}^{\vphantom{#3}}#2^{#3}}

\maketitle

\begin{abstract}
The concept of a {\em weight} on a planar algebra was introduced in 
\cite{DGG}. In this article we give an alternate characterization of weights on a planar algebra in terms of `weight functions' on the vertices of the principal graphs.
\comments{Using the weight functions, we show that 
if two bifinite bimodules have their bimodule planar algebras satisfying `trivial perturbation class', then so does their Connes fusion.}
Using this characterization we show that the property of bifinite bimodules of having a `trivial perturbation class' is closed under 
Connes fusion. 
We give a direct and constructive method of perturbing a bifinite bimodule by a positive weight in such a way that the bimodule planar algebra of the perturbed bimodule is isomorphic to the perturbation of the one associated to the initial bimodule by the given weight. 
\end{abstract}
\section{Introduction}
Starting from a $1$-cell in a pivotal $2$-category, a  purely algebraic construction of Jones' planar algebra was given by the second named author in \cite{Gho08} where the description of the action of tangles were given in terms of {\em graphical calculus of morphisms}, analogous to the ones used in \cite{Kas}.
In the operator algebra context, a nice prototype for this is the $2$-category of bifinite bimodules over $II_1$-factors; this was analyzed in great detail in \cite{DGG} where the concepts of {\em weight} on a planar algebra and {\em perturbation of a planar algebra by a weight} were introduced in order to illustrate the exact dependence of the planar algebra on the pivotal structure of the corresponding $2$-category.
Further, in the same paper, a correspondence was obtained between
a class of planar algebras slightly more general than subfactor planar algebras, referred as {\em bimodule planar algebra} (BMPA) on one hand and bifinite bimodules on the other, thereby providing a generalization of Jones' Theorem.  
\comments{
Further, in the same paper, a {\em bimodule planar algebra} (BMPA) was associated to each bifinite bimodule in such a way that extremality of the bimodule corresponds to sphericality of the planar algebra; conversely, every BMPA does come from a bifinite bimodule.}
Moreover, the {\em perturbation class} of a BMPA contains a unique spherical unimodular BMPA which can also be characterized by the property of having the lowest value of the index in the entire perturbation class.

This article grew out of the following questions raised in \cite{DGG} on perturbations of planar algebras.
\begin{enumerate}
\item Is there a way of characterizing the positive weights of a BMPA in terms of its principal graph?
\item 
Does the BMPA of the Connes fusion of two bimodules whose associated BMPAs have {\em trivial perturbation class} (TPC) (that is, every member of the perturbation class is spherical), have TPC?
In the subfactor context, this same as asking whether composition of TPC subfactors has TPC.
\item Given a bifinite bimodule, is there a `direct' way of perturbing it to another one by a weight in such a way that the BMPA associated to the perturbed bimodule is isomorphic to the perturbation of the BMPA associated to the initial bimodule by the weight?
Note that one can associate a bifinite bimodule to the perturbation of a BMPA associated to a bifinite bimodule by a weight using \cite[Theorem 5.13]{DGG} but it is unclear how the relation between the two bimodules depends on the weight.
\end{enumerate}
All the questions have been answered affirmatively in this paper. 
\vspace*{1em}

We now briefly describe the organization of this paper.
In Section \ref{prelim}, we set up some notations and give a quick recollection of various definitions and some standard facts on weights, perturbations and BMPAs.

Section \ref{wtfun} gives an alternate characterization of a positive weight on a BMPA in terms of {\em weight function} defined on the principal graphs which satisfies the {\em tensor homomorphism property}.
We further show that this weight function essentially depends on its value over the even vertices of any one of the principal graphs upto a positive scalar multiple of the weight function restricted to the odd vertices.
At the end of the section, we discuss a couple of examples of weight functions.

In Section \ref{pertcl}, using the weight functions, we affirmatively answer Question 2 above.
As a consequence, (i) Bisch-Haagerup planar algebras (see \cite{BDG09}) and (ii) $n$-cabling of any irreducible bimodule planar algebra for $n \geq 2$ always have TPC.

The final section starts with some notations and facts about imprimitivity bimodules and proofs of some of their key properties.
Using these bimodules and properties listed, we describe a method of perturbing a bifinite bimodule by a positive weight and prove that the BMPA of the new bimodule is indeed isomorphic to the pertubation of the BMPA corresponding to the old one by the given weight.
Moreover, if the initial bimodule is over hyperfinite $II_1$-factors, so is the perturbed bimodule, which is not the case for the bimodule that we get from \cite[Theorem 5.13]{DGG}.
\section{Preliminaries}\label{prelim}
\begin{figure}[h]
\psfrag{k}{$k$}
\psfrag{2k}{$2k$}
\psfrag{e}{$\vlon$}
\psfrag{-e}{$-\vlon$}
\psfrag{M}{$ M_{\vlon k} =$}
\psfrag{m}{$(\vlon k,\vlon k)\rightarrow \vlon k$}
\psfrag{1ek}{$1_{\vlon k} =$}
\psfrag{1}{$\emptyset \rightarrow \vlon k$}
\psfrag{RI}{$RI_{\vlon k} =$}
\psfrag{ri}{$\vlon k \rightarrow \vlon (k+1)$}
\psfrag{LI}{$LI_{\vlon k}=$}
\psfrag{li}{$\vlon k \rightarrow -\vlon (k+1)$}
\psfrag{Id}{$I_{\vlon k}=$}
\psfrag{id}{$\vlon k \rightarrow \vlon k$}
\psfrag{E}{$E_{\vlon (k+1)}=$}
\psfrag{jp}{$\emptyset \rightarrow \vlon (k+2)$}
\psfrag{RE}{$RE_{\vlon (k+1)} =$}
\psfrag{re}{$\vlon(k+1) \rightarrow \vlon k$}
\psfrag{LE}{$LE_{\vlon (k+1)} =$}
\psfrag{le}{$\vlon(k+1) \rightarrow -\vlon k$}
\includegraphics[scale=0.21]{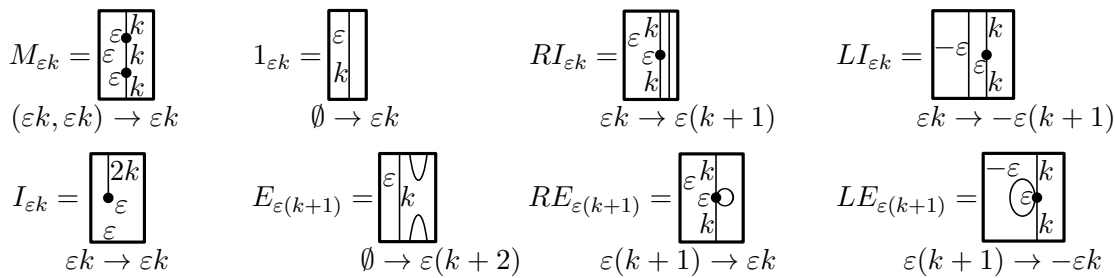}\caption{Generating tangles.}\label{tangles}
\end{figure}

In this section, we will recall certain basic facts about planar algebras which will be used in the forthcoming sections.
For planar algebras, we will stick to the notations and terminology set up in Section $2.1$ of \cite{DGG}; however, we would like to briefly recall the key points in order to make this article self-contained.
\begin{enumerate}
\item We will consider the natural binary operation on $\{-,+\}$ given by $++ := +$, $+-:=-$, $-+:=-$ and $--:=+$.
\item We will denote the set of all possible colors of discs in tangles by $Col:= \left\{ \vlon k: \vlon \in \{+,-\},\, k \in \mbb{N}_0\right\}$ where $\N_0 := \N \cup \{ 0 \}$.
\item In a tangle, we will replace (isotopically) parallel strings by a single strand labelled by the number of strings, and an internal disc with color $\vlon k$ will be replaced by a bold dot with the sign $\vlon$ placed at the angle corresponding to the distinguished boundary components of the disc.
For example,\\
\psfrag{replace}{will be replaced by}
\psfrag{2}{$2$}
\psfrag{4}{$4$}
\psfrag{e}{$\vlon$}
\includegraphics[scale=0.21]{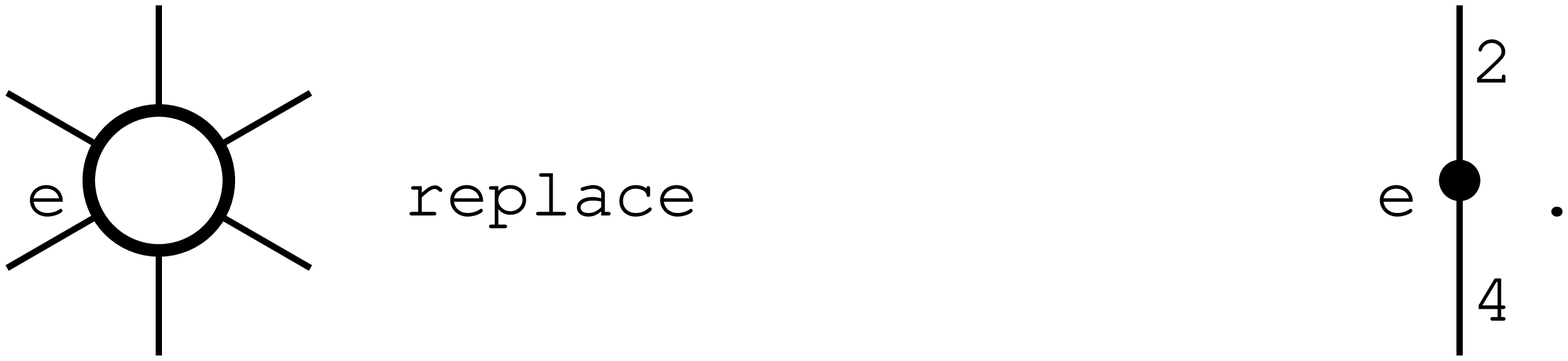}
\end{enumerate}
We set up some notation for a set of `generating tangles' in Figure \ref{tangles}.

If $P$ is a planar algebra, then by a $P$-{\em labelled} (resp., {\em semi-labelled}) tangle, we mean a tangle, all (resp., some) of whose internal discs are labelled by elements of $P$ in such a way that an internal disc of color $\vlon k$ is labelled by an element of $P_{\vlon k}$.
For simplicity, we will replace a $P$-labelled internal disc by a bold dot as before, with the label being placed at the angle corresponding to the distinguished boundary component of the disc.
Also, ${\mcal T}_{\vlon k}$ (resp., ${\mcal T}_{\vlon k} (P)$) will denote the set of tangles (resp., $P$-labelled tangles) which has $\vlon k$ as the color of the external disc; ${\mcal P}_{\vlon k} (P)$ will be the vector space with ${\mcal T}_{\vlon k}$ as a basis.
We continue to denote the linear map induced by the action of $P$ by $P : {\mcal P}_{\vlon k} (P) \rightarrow P_{\vlon k}$.
\vspace*{2mm}

We now recall the definitions of {\em weights} and {\em perturbations} of planar algebras, introduced in \cite{DGG}, and a few general facts.
\begin{defn}\label{weight}
Let $P$ be a planar algebra. An invertible element $z \in P_{+1} $ is said to be a weight of $P$ if $z_{\vlon k} \in \mcal{Z}(P_{\vlon k})$ for all $\vlon k \in Col$, where
\begin{equation}\label{zek}
\psfrag{z}{$z$}
\psfrag{zinv}{$z^{-1}$}
\psfrag{kstrings}{$\cdots k$\text{ strings}}
\psfrag{+}{$+$}
\psfrag{-}{$-$}
z_{+k} := P_{
\includegraphics[scale=0.3]{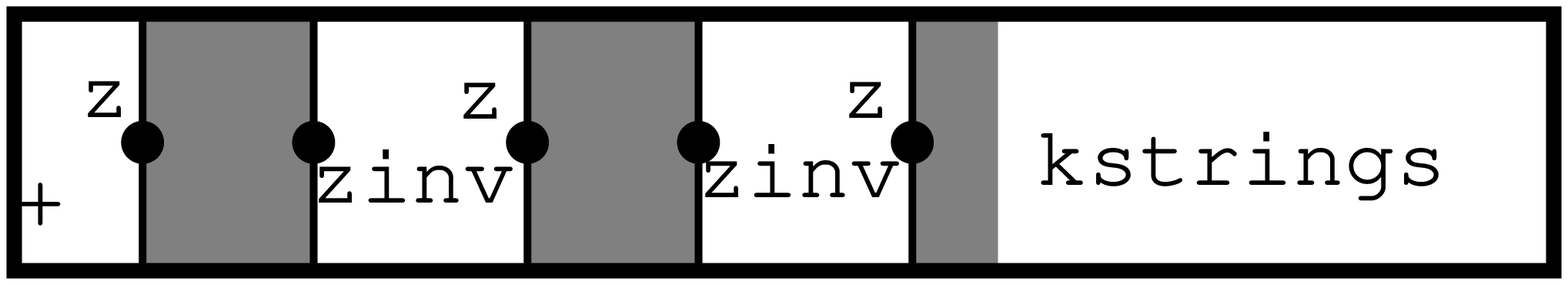}
}
\text{ and } z_{-k} := P_{
\includegraphics[scale=0.3]{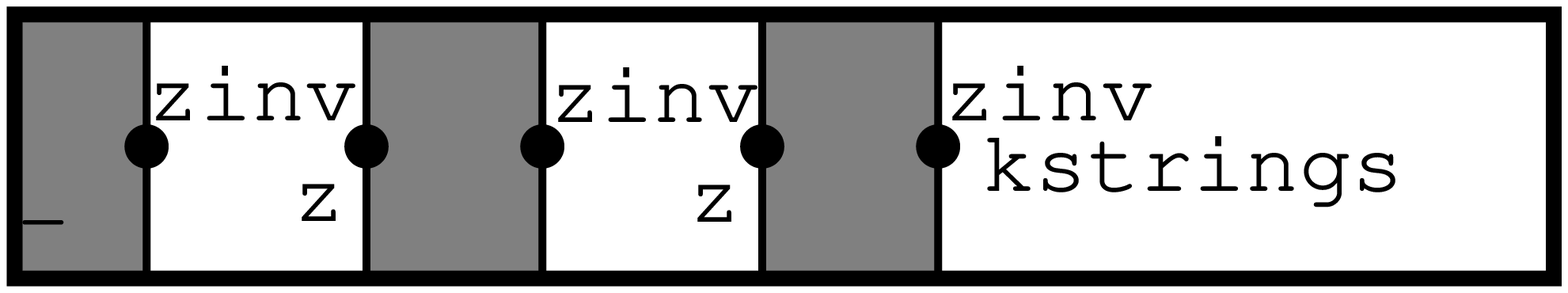}
}.
\end{equation}
\end{defn}
Given a weight $z$ of a planar algebra $P$ and an invertible decomposition $z = ab$ for $a,b$ invertible in $P_{+1}$, we now construct a new planar algebra $P^{(a,b)}$ as follows:

\noindent (i) {\em Vector spaces}: $P^{(a,b)}_{\vlon k} : = P_{\vlon k}$ for all $ \vlon k \in Col$.

\begin{figure}[h]
\psfrag{a}{$a$}
\psfrag{b}{$b$}
\psfrag{ainv}{$a^{-1}$}
\psfrag{binv}{$b^{-1}$}
\psfrag{by}{by}
\psfrag{;}{;}
\includegraphics[scale=0.15]{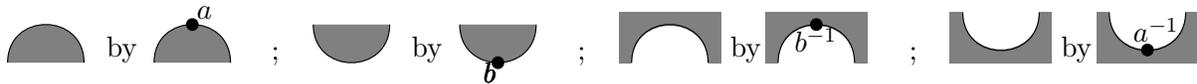}
\caption{Perturbing a planar algebra}
\label{max-min}
\end{figure}

\noindent (ii) {\em Actions of tangles}: Let $T$ be a tangle and $\hat{T}$ be a standard form representative (see \cite[$\S$4]{Gho08}) of the isotopy class of $T$.
We replace each local maximum and minimum appearing in $\hat T$ according to the prescription given by Figure \ref{max-min} and call the resulting semi-labelled tangular diagram $\hat{T}^{(a,b)}$.
Define $P^{(a,b)}_T := P_{\hat{T}^{(a,b)}}$.\\
It can be checked that $P^{(a,b)}_T$ is indeed well-defined.

Note that $P^{(a,b)}$ has same filtered algebra structure as $P$ whereas the action of Jones projection tangles and conditional expectation tangles differ.
We will refer $P^{(a,b)}$ as the {\em perturbation of $P$ by the decomposition $z = ab$ of the weight $z$}.
\begin{rem}
For an invertible decomposition $z = a b$ of a weight $z$ of a planar algebra $P$ and any $\lambda \in \C \setminus \{ 0\}$, $P^{(a,b)} = P^{(\lambda a, \lambda^{-1}b)}$. Further, the planar algebras $P^{(a, b)}, P^{(b, a)}, P^{(z, 1)},$ and $ P^{(1, z)}$ are all isomorphic.
Hence, up to isomorphism, the perturbation of $P$ only depends upon the weight $z$.
\end{rem}

A trivial example of a weight of a planar algebra $P$ is a non-zero scalar $\lambda \in \C $.
By the above remark, $P^{(\lambda, \mu)} = P^{(\lambda \mu, 1)} = P^{(1, \lambda \mu)} = P^{(\mu, \lambda)}$ for all non-zero scalars $\lambda$ and $\mu$.
We will usually refer to such perturbations as {\em scalar perturbations}.
A perturbation of a planar algebra with modulus need not have modulus (except for perturbations by scalar weights).
However, if the planar algebra is connected then so are its perturbations, but the moduli of the perturbations might vary.
For instance, for a planar algebra $P$ with modulus $(\delta_-, \delta_+)$, the scalar perturbation $P^{(\lambda, 1)}$ has modulus $( \lambda^{-1} \delta_-,\, \lambda \delta_+ )$.
\begin{defn}
The normalized planar algebra associated to a planar algebra $P$ with modulus $(\delta_-, \delta_+)$ is its scalar perturbation by the weight $ \sqrt{\frac{\delta_-}{\delta_+}}$.
\end{defn}
\noindent Note that the normalization of  $P$ is a unimodular planar
algebra.  Although scalar perturbations change the modulus, the index
however remains the same.
\begin{defn}
A connected planar algebra $P$ is said to be spherical if the actions of $0$-tangles in its normalization are invariant under spherical isotopy.
\end{defn}
Note that this property is equivalent to demanding that the normalized left and right picture traces on $P_{+1}$ are identical. In general, a perturbation of a $\ast$-planar algebra need not be a $\ast$-planar algebra.
However, for certain specific weights, the perturbations also turn out to be $\ast$-planar algebras.
For instance, if $P$ is a $\ast$- (resp., positive) planar algebra, it is routine to verify that a perturbation of the type $P^{(a, \lambda a^*)}$ for non-zero real (resp., positive) scalar $\lambda$ becomes a $\ast$- (resp., positive) planar algebra with $\ast$-structure coming from the original one.
\vspace*{2mm}

We now state the `Jones' theorem for bimodules' which provides an important link between bifinite bimodules over $II_1$-factors and planar algebras.
A finite dimensional, connected, positive $C^*$-planar algebra will be referred as a {\em bimodule planar algebra}.
Let $\lsub{A_+}{\mathcal H}_{A_-}$ be a finite index bimodule for $II_1$ factors $A_\pm$.
Before stating the theorem, we set up some more notations:\\
${\mathcal H}_{+ 0} := L^2 (A_+)$, ${\mathcal H}_{- 0} : = L^2 (A_-)$ and for $k\geq 1$, ${\mathcal H}_{\vlon k}$ is the tensor product (over $A_+$ or $A_-$) of $k$-many modules $\mathcal H$ and $\overline{\mathcal H}$ alternately with $\mathcal H$ (resp., $\overline{\mathcal H}$) being the left-most module if $\vlon=+$ (resp., $-$).
\begin{thm}\label{bimod-pa-theorem} (see \cite{DGG})
(i) $P$ defined by $P_{\vlon k} := \lsub{A_\vlon}{\mathcal L}_{A_\eta} ({\mathcal H}_{\vlon k})$ where $\eta = (-)^k \vlon$ , has a unique bimodule planar algebra structure with $\ast$-structure coming from the usual adjoints of intertwiners, satisfying:\\
(a) action of multiplication tangles matches with the composition of operators in the intertwiner spaces,\\
(b) $P_{RI_{\vlon k}} \! = \! \left\{ \! \! \! \!
\begin{tabular}{ll}
$P_{\vlon k}\ni T \mapsto T \underset{A_+}{\otimes} id_{\mathcal H} \in P_{\vlon (k+1)}$, & if either $\vlon = +$ and $k$ is even, or $\vlon = -$ and $k$ is odd,\\
$P_{\vlon k}\ni T \mapsto T \underset{A_-}{\otimes} id_{\overline{\mathcal H}} \in P_{\vlon (k+1)}$, & otherwise,
\end{tabular}
\right.$\\
(c) $P_{LI_{\vlon k}} \! = \! \left\{ \! \! \! \!
\begin{tabular}{ll}
$P_{\vlon k}\ni T \mapsto id_{\overline{\mathcal H}} \underset{A_+}{\otimes} T \in P_{-\vlon (k+1)}$, & if $\vlon = +$,\\
$P_{\vlon k}\ni T \mapsto id_{\mathcal H} \underset{A_-}{\otimes} T \in P_{-\vlon (k+1)}$, & if $\vlon = -$, and
\end{tabular}
\right.$\\
(d) $P_{E_{+1}}$ (resp., $P_{E_{-1}}$) is given by 
\begin{align*}
& {\mathcal H}^o \underset{A_-}{\otimes} {\overline{\mathcal H}}^o \ni \xi \underset{A_-}{\otimes} \overline{\eta} \overset{A_+ \text{-} A_+}{\longmapsto} \underset{j}{\sum} \; \lsub{A_+}{\langle \eta , \xi \rangle} \left( \eta_j \underset{A_-}{\otimes} \overline{\eta}_j \right) \in {\mathcal H}^o \underset{A_-}{\otimes} {\overline{\mathcal H}}^o\\
\text{(resp., } & {\overline{\mathcal H}}^o \underset{A_+}{\otimes} {\mathcal H}^o \ni \overline{\eta} \underset{A_+}{\otimes} \xi \overset{A_- \text{-} A_-}{\longmapsto} \underset{i}{\sum} {\langle \eta , \xi \rangle}_{A_-} \left( \overline{\xi}_i \underset{A_+}{\otimes} \xi_i \right) \in {\overline{\mathcal H}}^o \underset{A_+}{\otimes} {\mathcal H}^o \text{ )}
\end{align*}
where $\{\xi_i\}_{i}$ (resp., $\{\eta_j\}_{j}$) is any basis for $\lsub{A_+}{\mathcal H}$ (resp., ${\mathcal H}_{A_-}$).
(We will refer $P$ (resp., normalization of $P$) as the {\em bimodule} (resp., {\em normalized} or {\em unimodular}) {\em planar algebra associated to $_A{\mathcal H}_B$}).

(ii) $P$ is spherical if and only if $\lsub{A_+}{\mathcal H}_{A_-}$ is extremal.

(iii) If $B := {\mathcal L}_{A_-} ({\mathcal H})$, then the normalized planar algebras associated to $\lsub{A_+}{\mathcal H}_{A_-}$ and $\lsub{A_+}{L^2 (B)}_B$ are isomorphic as $\ast$-planar algebras.

(iv) Given any bimodule planar algebra $P$, there exists a bifinite bimodule whose associated bimodule planar algebra is isomorphic to $P$.
\end{thm}
\section{An alternate characterization of weights}\label{wtfun}
In this section, we will first find a necessary and sufficient condition for existence of a positive weight on a bimodule planar algebra in terms of its principal graphs; this answers a question in \cite{DGG}.
This was suggested by Stefaan Vaes. We further show that this condition mainly depends on the even vertices of any one of the principal graphs.
In the end, we discuss some examples.

Let $P$ be the bimodule planar algebra associated to the bifinite bimodule ${_{A_+}} {\mcal H}_{A_-}$ for $II_1$-factors $A_{\pm}$ as in Theorem \ref{bimod-pa-theorem} and $z$ be a positive weight on $P$.
We will refer an $A_\vlon$-$A_\eta$ bimodule as $\vlon$-$\eta$ bimodule for $\vlon , \eta \in \{- , +\}$.
Set
\[
V_{\vlon , \eta} := \{ \text{isomorphism classes of irreducible } \eta \text{-} \vlon \text { sub-bimodule of } {\mcal H}_{\eta k} : k \in (2\N_0 + \delta_{\vlon \neq \eta}) \}.
\]
Note that for each $v \in V_{\vlon , \eta}$, there exists a $k \in (2\N_0 + \delta_{\vlon \neq \eta})$ and $p \in {\mathscr P}_{\text{min}} (P_{\eta k})$ such that $v = [\text{Range } p]$.
Now, since $z_{\eta k}$ is central, for each $p \in {\mathscr P}_{\text{min}} (P_{\eta k})$, there exists $w_p > 0$ such that $z_{\eta k} p = w_p p$.
Fix $\vlon, \eta \in \{- , +\}$.
Define the function 
\[
V_{\vlon , \eta} \ni v \os{w}{\longmapsto} w_v := w_p \in (0,\infty) \text{ where $p \in {\mathscr P}_{\text{min}} (P_{\eta k})$ such that $v = [\text{Range } p]$.}
\]
We need to check that $w$ is well-defined. Suppose $q \in {\mathscr P}_{\text{min}} (P_{\eta k+2l})$ for $l \geq 0$ such that $v = [\text{Range } q]$.
It is enough to check $w_p = w_q$.
Now, from minimality of $p$ and the action of Jones projection tangle in Theorem \ref{bimod-pa-theorem}, we may conclude:

(i) $p_1 := \delta^{-1}_{\vlon} p P_{E_{\eta (k+1)}} \in {\mathscr P}_{\text{min}} (P_{\eta (k+2)})$,

(ii) $\text{Range } p_1 \cong \left(\text{Range } p \us{\vlon}{\otimes} L^2 (A_\vlon) \right) \cong \text{Range } p \cong \text{Range } q$ as $\eta$-$\vlon$ bimodules where $\us{\vlon}{\otimes}$ denotes fusion over $A_\vlon$ and

(iii) $w_{p_1} = w_p$ since $w_{p_1} p_1 = z_{\eta (k+2)} p_1 = \delta^{-1}_{\vlon} P_{RI^2_{\eta k} (z_{\eta k} p)} P_{E_{\eta (k+1)}} = w_p p_1$ where $RI^m_{\eta k} : \eta k \rightarrow \eta (k+m)$ is the tangle obtained from $RI_{\eta k}$ by replacing the single string on the right of the internal disc by $m$ parallel strings.\\
Iterating the above, we get $p_{l} := \delta^{-l}_{\vlon} P_{RI^{2l}_{\eta k}} (p) P_{E_{\eta (k+1)}} P_{E_{\eta (k+3)}} \cdots P_{E_{\eta (k+2l-1)}} \in {\mathscr P}_{\text{min}} (P_{\eta (k+2l)})$ satisfying $w_{p_l} = w_p$ and $[\text{Range } p_l] = [\text{Range } q] = v$.
Thus, $q$ and $p_l$ are MvN-equivalent in $P_{\eta (k+2l)}$ and hence, $w_q = w_{p_l} = w_p$.

Further, note that if $v_1 \in V_{\eta , \nu}$, $v_2 \in V_{\vlon , \eta}$ and $v_3 \in V_{\vlon , \nu}$ such that $v_3 \leq v_1 \us{\eta}{\otimes} v_2$, then using the action of the left and the right inclusion tangles, we get
\[
w_{v_3} p_{3} = z_{\nu (k+l)} p_3 = \left( z_{\nu k} \us{\eta}{\otimes} z_{\eta l} \right) (p_1 \us{\eta}{\otimes} p_2) p_3 = w_{v_1} w_{v_2} (p_1 \us{\eta}{\otimes} p_2) p_3 = w_{v_1} w_{v_2} p_3
\]
where $p_1 \in {\mathscr P}_{\text{min}} (P_{\nu k})$, $p_2 \in {\mathscr P}_{\text{min}} (P_{\eta l})$ and $(p_1 \us{\eta}{\otimes} p_2) \geq p_3 \in {\mathscr P}_{\text{min}} (P_{\nu (k+l)})$ such that $v_i = [\text{Range } p_i]$ for all $i = 1, 2, 3$. This motivates the following definition.
\begin{defn}
$w$ is said to be a weight function for a semisimple bicategory ${\mcal B}$ if for all $A, B \in {\mcal B}_0$, there exists a map $w : V_{A,B} := \{ \text{isomorphism classes of simple objects in }  {\mcal B}_{A,B} \} \rightarrow (0, \infty)$ satisfying
\begin{equation}\label{thp}
w_{v_3} = w_{v_1} w_{v_2} \text{ whenever } v_1 \otimes v_2 \text{ contains } v_3
\end{equation}
for all $v_1 \in V_{B , C}$, $v_2 \in V_{A , B}$ and $v_3 \in V_{A , C}$.
\end{defn}
Equation \ref{thp} will be referred as {\em tensor homomorphism property}. We establish the connection between weights and weight functions in the following lemma.
\begin{lem}\label{wtchar}
There is a one-to-one correspondence between positive weights on the planar algebra $P$ associated to a bifinite bimodule ${_{A_+}} {\mcal H}_{A_-}$ for $II_1$-factors $A_{\pm}$, and weight functions on the bicategory of bimodules generated by ${_{A_+}} {\mcal H}_{A_-}$.
\end{lem}
\begin{pf}
We have already associated a weight function to every positive weight on $P$.
For the converse, let $w$ be a weight function for the bicategory.
Set $z_{\vlon k} := \us{p \in {\mathscr P}_{\text{min}} ({\mathcal Z} (P_{\vlon k}) ) }{\sum} w_{v_p} p$ where $v_p = [ \text{Range } p_0 ]$ for any minimal projection under $p$.
Clearly, $z_{\vlon k}$ is a central, invertible, positive element in $P_{\vlon k}$.
So, it remains to establish Equation \ref{zek} with $z$ replaced by $z_{+1}$.
\comments{
\begin{equation}
\psfrag{z}{$z$}
\psfrag{zinv}{$z^{-1}$}
\psfrag{kstrings}{$\cdots k$\text{ strings}}
\psfrag{+}{$+$}
\psfrag{-}{$-$}
z_{+k} = P_{
\includegraphics[scale=0.3]{figures/newchar/z+k.eps}
}\text{ and } z_{-k} = P_{
\includegraphics[scale=0.3]{figures/newchar/z-k.eps}
}
\end{equation}
where $z = z_{+1}$.
}
For this, note that the tensor homomorphism property of any weight $w$ requires it to assign $1$ to the unit objects, and thereby, $w_{\oline{v}} = w^{-1}_v$ for all $v \in V_{\vlon , \eta}$, $\vlon , \eta \in \{- , +\}$.
We will use induction on $k$.\\
Suppose $k=1$.
We need to show $z_{-1} = P_{R^1_{+1}} (z^{-1})$. Observe that right side of this equation is also central since $P_{R^1_{+1}}$ is an anti-algebra isomorphism. Also, if $v = [\text{Range } p]$ for $p \in {\mathscr P}_{\text{min}} (P_{+1})$, then $\oline{v} = [\text{Range } P_{R^1_{+1}} (p)]$.
Thus, $P_{R^1_{+1}} (z^{-1}) P_{R^1_{+1}} (p) = P_{R^1_{+1}} (w^{-1}_v p) = w_{\oline{v}} P_{R^1_{+1}} (p) = z_{-1} P_{R^1_{+1}} (p)$ which gives the required equation.\\
Suppose the equations hold for $k = l$ and we will check that they also hold for $k = l+1$.
By the induction assumption and action of the left and the right inclusion tangles, the right hand side of the first equation in Equation \ref{zek} can be expressed as
\[
z_{+l} \us{\vlon}{\otimes} z_{\vlon 1}
= \us{i,j}{\sum} w_{v_{p_i}} w_{v_{q_j}} \left(p_i \us{\vlon}{\otimes} q_j \right)
= \us{i,j}{\sum} \us{s(i,j)}{\sum} w_{[\text{Range } r_{s(i,j)}]} r_{s(i,j)} = z_{+(l+1)}
\]
where $\vlon = (-)^l$, $\{p_i\}_i = {\mathscr P}_{\text{min}} ({\mcal Z} (P_{+l}) )$, $\{q_j\}_j = {\mathscr P}_{\text{min}} ({\mcal Z} (P_{\vlon 1}) )$ and $\{ r_{s(i,j)} \}_{s(i,j)}$ is an orthogonal decomposition of the projection $\left(p_i \us{\vlon}{\otimes} q_j \right)$ into minimal ones. 
The equation on the right side can also be deduced using similar arguments.
\end{pf}

One can also consider weight functions for semisimple tensor categories instead of bicategories.
It is natural to ask whether the restriction of a weight function for the bicategory generated by ${_{A_+}} {\mcal H}_{A_-}$ to the tensor category of $A_+$-$A_+$ bimodules, contain all information about $w$.
We answer this question affirmatively in the next proposition.
\begin{prop}\label{wtext}
For $\vlon = \pm$, any weight function $w$ for the tensor category of $A_\vlon$-$A_\vlon$ bimodules generated by ${_{A_+}} {\mcal H}_{A_-}$, can be extended to a weight function for the bicategory associated to  ${_{A_+}} {\mcal H}_{A_-}$.
\end{prop}
\begin{pf}
Without loss of generality, let $\vlon = +$.
Set $w_{+ 2l} := \us{p \in {\mathscr P}_{\text{min}} ({\mathcal Z} (P_{+ 2l}) ) }{\sum} w_{v_p} p$ where $v_p = [ \text{Range } p_0 ]$ for any minimal projection under $p$.
We first deduce the following useful equations involving $w_{+2}$:

(i) $P_{R^2_{+2}} (w_{+2}) = w^{-1}_{+2}$,

(ii) $w_{+2l} = P_{
\psfrag{w}{$w_{+2}$}
\psfrag{l}{$\cdots 2l$ strings}
\psfrag{2}{$2$}
\psfrag{+}{$+$}
\includegraphics[scale=0.25]{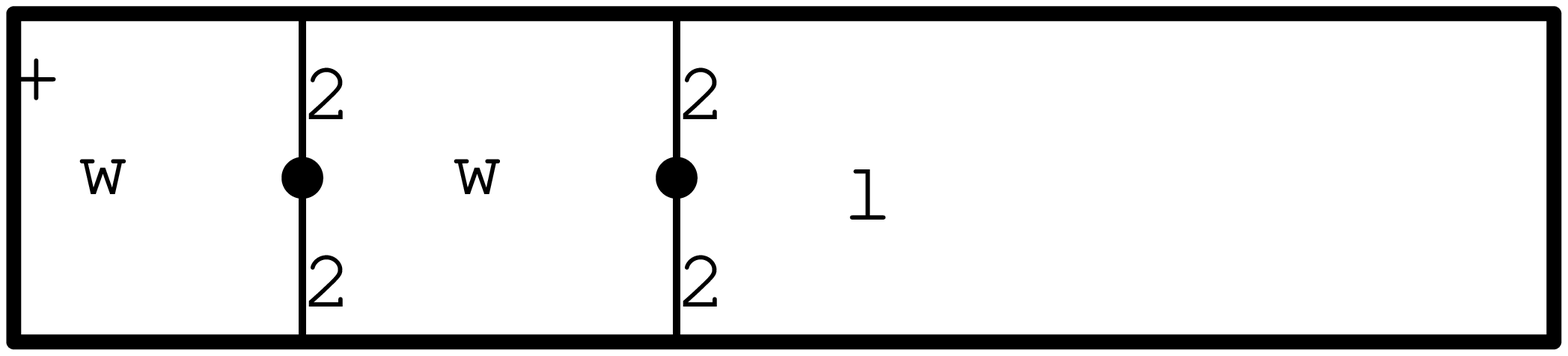}
}$ ,

(iii) $P_{
\psfrag{w}{$w_{+2}$}
\psfrag{+}{$+$}
\includegraphics[scale=0.2]{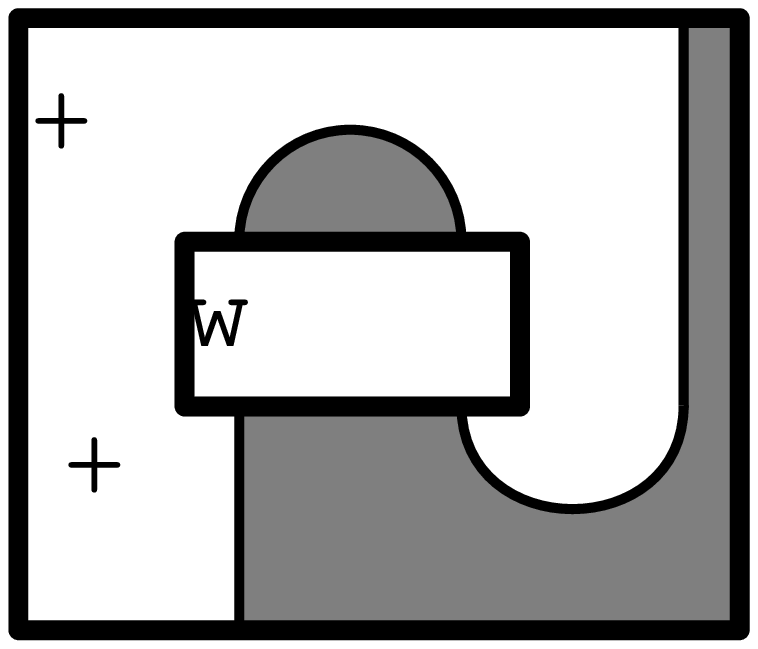}
} 
= 1_{P_{+1}} = 
P_{
\psfrag{w}{$w_{+2}$}
\psfrag{+}{$+$}
\includegraphics[scale=0.2]{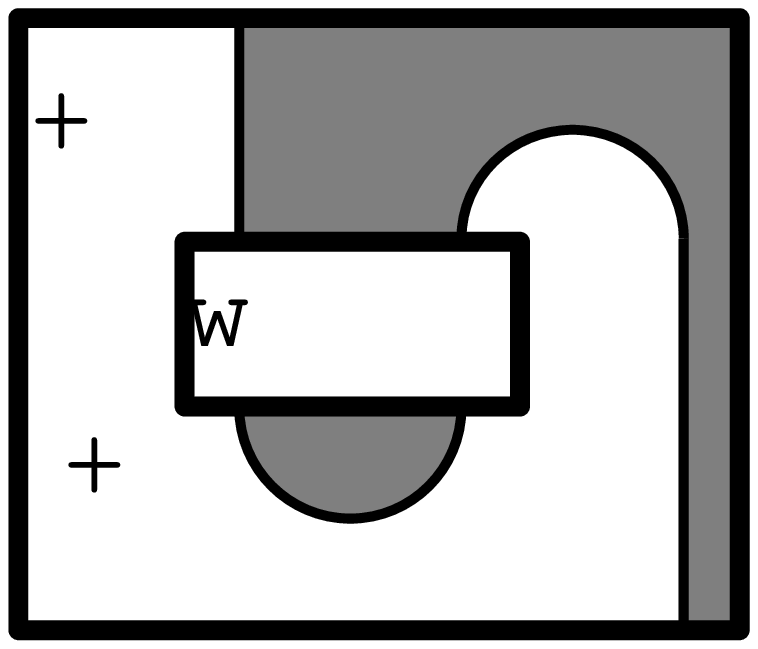}
}$ ,

(iv) {\em Exchange relations:}
$P_{
\psfrag{w}{$w_{+2}$}
\psfrag{+}{$+$}
\includegraphics[scale=0.2]{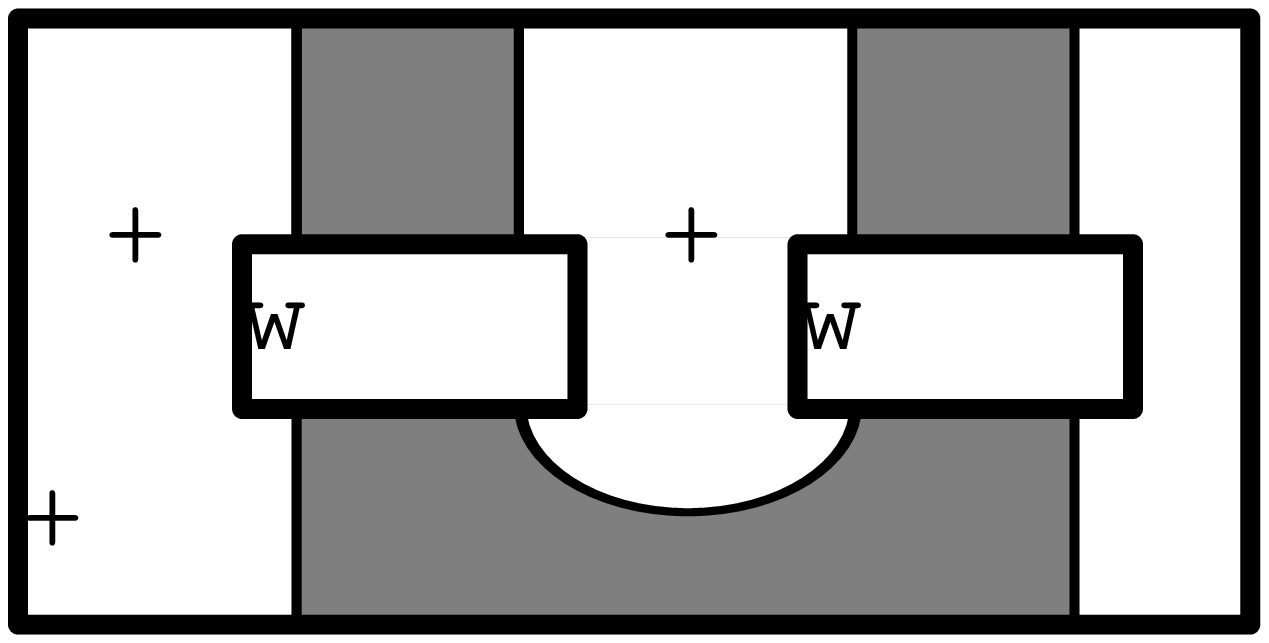}
}
= P_{
\psfrag{w}{$w_{+2}$}
\psfrag{+}{$+$}
\includegraphics[scale=0.2]{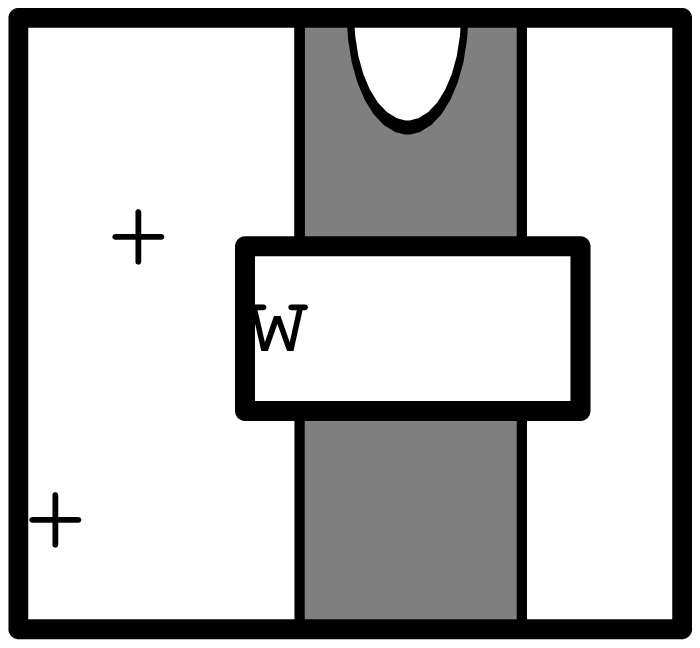}
}$ and 
$P_{
\psfrag{w}{$w_{+2}$}
\psfrag{+}{$+$}
\includegraphics[scale=0.2]{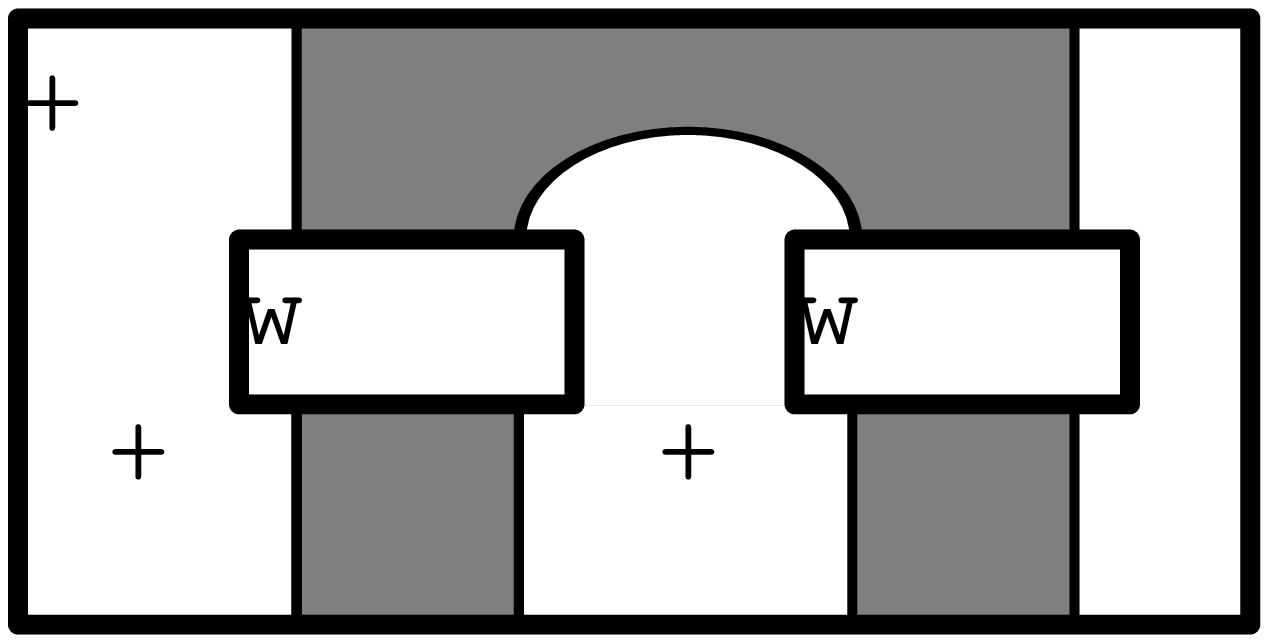}
} = P_{
\psfrag{w}{$w_{+2}$}
\psfrag{+}{$+$}
\includegraphics[scale=0.2]{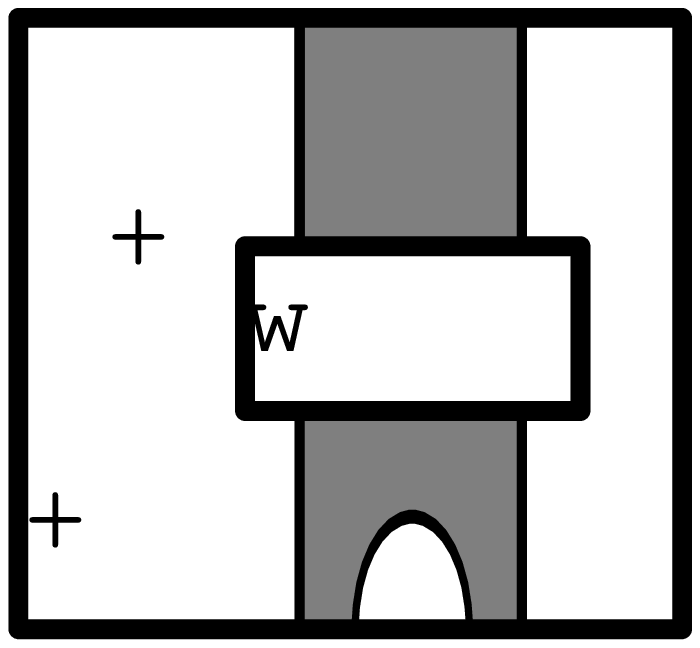}
}$ .

From the action of rotation tangle, one can show $\text{Range } P_{R^2_{+2}} (p) \cong \oline{\text{Range } p}$ as $A_+$-$A_+$ bimodules for all $p \in {\mscr P}(P_{+2})$; this along with the fact that $w_{\oline{v}} = w^{-1}_v$ for all $v \in V_{+,+}$, implies part (i).

Part (ii) easily follows from the tensor homomorphism property of $w$.

For part (iii), it is enough to show that $w_{+2} P_{E_{+2}} = P_{E_{+2}}$. Note that $\delta^{-1}_+ P_{E_{+2}}$ is a minimal projection of $P_{+2}$ whose range is $L^2 (A_+)$. Hence, we have (iii).

Part (iv) can derived from the equation $[w_{+4} , P_{E_{+2}} P_{E_{+3}}]=0$ ($w_{+4}$ being central), (ii) and (iii).
\vspace{1em}

Consider the positive element $z = P_{RE_{+2}} (z_{+2}) \in P_{+1}$ and set $\lambda := P_{LE_{+1}} (z) \in \C \simeq P_{-0}$.
Note that by equations (i) and (iv), we have $z^{-1} = \lambda^{-1} P_{RE_{+2}} (w^{-1}_{+2})$.
By Lemma \ref{wtchar}, it is enough to show that $z$ is a positive weight for $P$ and $z_{+2l} = w_{+2l}$ for all $l \in \N$ where $z_{\vlon k}$'s denote the usual elements obtained from $z$.
So,
\[
z_{+2} = \lambda^{-1} P_{
\psfrag{w}{$w_{+2}$}
\psfrag{+}{$+$}
\includegraphics[scale=0.2]{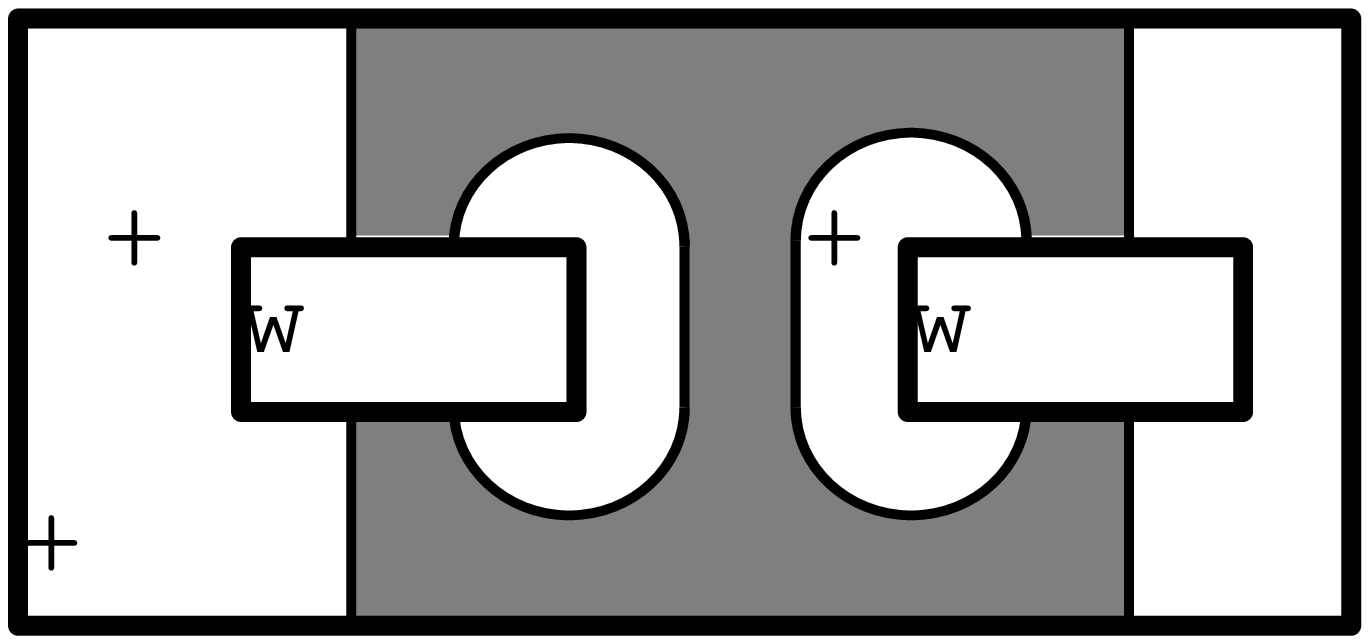}
} = \lambda^{-1} P_{
\psfrag{w}{$w_{+2}$}
\psfrag{+}{$+$}
\includegraphics[scale=0.2]{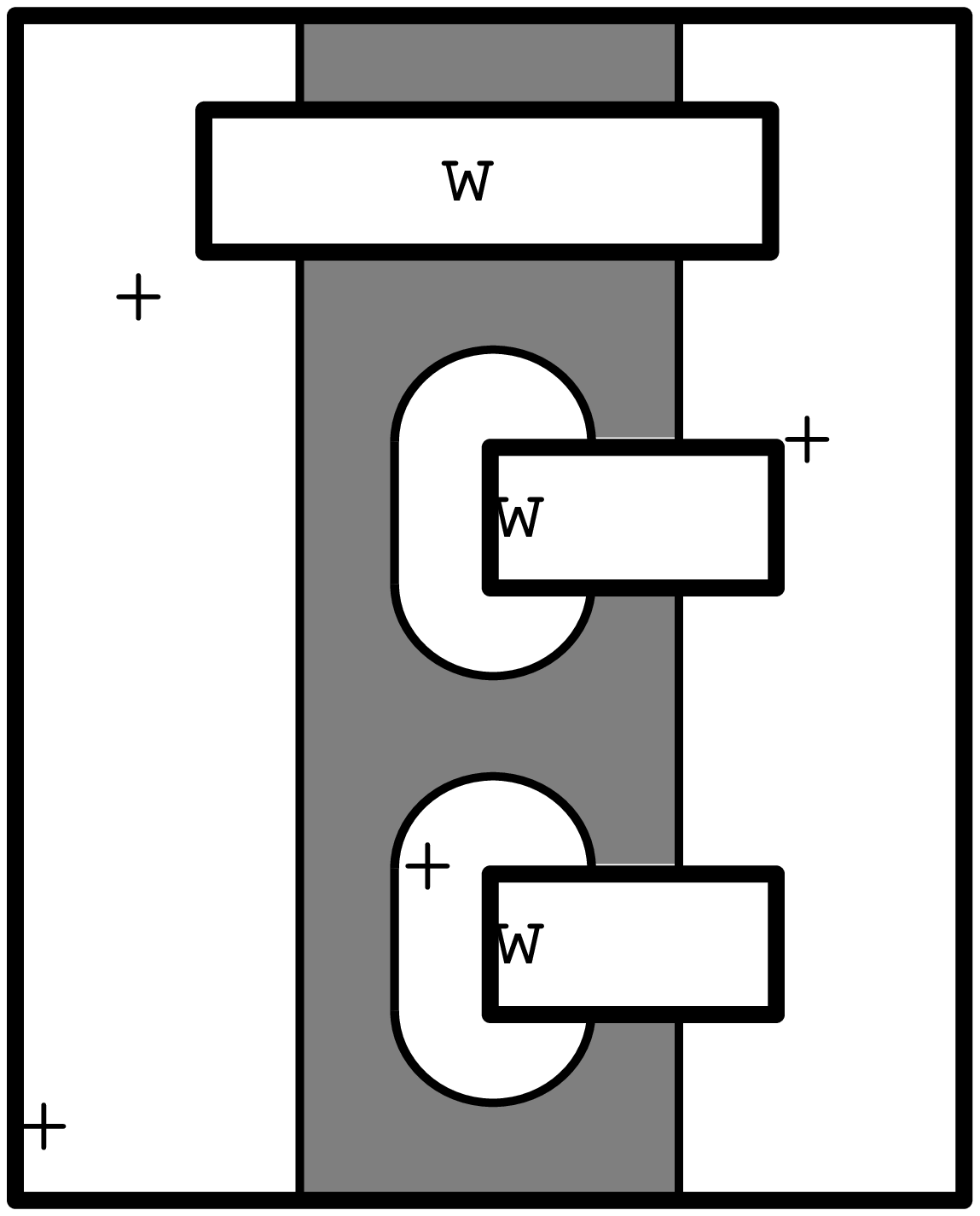}
} = w_{+2}
\]
where the second equality is obtained applying the first exchange relation in (iv) and the third comes from part (i) and the expression for $z^{-1}$.
Hence, by equation (ii), $z_{+2l} = w_{+2l}$ for $l \in \N$.

It remains to be shown that $z_{\vlon k}$'s are central.
The case $\vlon k = +2l$ is already done. Suppose $\vlon k = -2l$. Since the action of tangles preserves composition, we get
\begin{align*}
z_{-2l} \; x
& = \lambda^{-1} P_{LE_{+(2l+1)} \circ RE_{+2(l+1)} } \left( w_{+2(l+1)} \right) \; x \\
& = \lambda^{-1} P_{LE_{+(2l+1)} \circ RE_{+2(l+1)} } \left( w_{+2(l+1)} \; P_{LI_{-(2l+1)} \circ RI_{-2l} } (x) \right)\\
& = \lambda^{-1} P_{LE_{+(2l+1)} \circ RE_{+2(l+1)} } \left(P_{LI_{-(2l+1)} \circ RI_{-2l} } (x) \; w_{+2(l+1)} \right)\\
& = \lambda^{-1} x \; P_{LE_{+(2l+1)} \circ RE_{+2(l+1)} } \left( w_{+2(l+1)} \right)
= x \; z_{-2l}
\end{align*}
for all $x \in P_{-2l}$.
Hence, $z_{-2l}$ is central and the remaining cases can be shown similarly.
\end{pf}
\begin{rem}\label{wtextunique}
From the above proof, observe that the extension of the weight function $w$ (in Proposition \ref{wtext}) is unique only upto a multiple of its restriction to the irreducible $A_+$-$A_-$ bimodules, by a positive scalar.
\end{rem}
\begin{example}
The first easy example of a (non-trivial) weight function appears in the case $V_{+,+} = \langle a \rab \cong \Z$ where every choice of $\lambda > 0$ gives a weight function $a^n \mapsto \lambda^n$.
In the subfactor context, these weight functions correspond to the non-extremal subfactors constructed by Jones in \cite{Jon83}, whose associated planar algebras were constructed explicitly using weights and perturbations in \cite{DGG}.
\end{example}
\begin{example}\label{qgpeg}
Another example of a weight function comes from the representation category of the {\em universal unitary compact quantum group} $A_u (F)$ for $F \in GL(n,\C)$ (introduced in \cite{VW} and studied in detail in \cite{Ban}).
The set of isomorphic classes of irreducible representations of $A_u (F)$ is indexed by the free monoid $\N \ast \N = \lab \alpha \rab \ast \lab \beta \rab$ with the involution $\ol{\alpha} = \beta$.
The tensor product structure is given by
\[
w_1 \otimes w_1 \cong \us{
\left\{
(a,b) \in \N \ast \N
\left|\begin{array}{l}
\text{there exists } c \in \N \ast \N \text{ such}\\
\text{that } w_1 = ac \text { and } w_2 = \ol{c} b
\end{array}\right.
\right\}
}{\bigoplus} ab 
\]
where $w_1 , w_2 \in \N \ast \N$.
Then, for each $\lambda > 0$, we can define a weight function by $\N \ast \N \ni w \mapsto \lambda^{m - n}$ where $m$ (resp., $n$) is the number of $\alpha$ (resp., $\beta$) in $w$. This example was pointed out by Stefaan Vaes.
\end{example}
\section{Trivial perturbation class}\label{pertcl}
We begin this section by recalling the definition of `trivial perturbation class' (introduced in \cite{DGG}) and then, prove that this property is preserved under composition of subfactors.
This answers a question in \cite{DGG}.
\begin{defn}
A bimodule planar algebra is said to have trivial perturbation class (TPC) if all its perturbations by positive weights, are spherical.
\end{defn}
By a subfactor or a bimodule having TPC, we will mean that its associated planar algebra has TPC.
\begin{lem}\label{tpceq}
If $P$ is the bimodule planar algebra associated to a bifinite bimodule ${_{A_+}} {\mcal H}_{A_-}$ for $II_1$-factors $A_\pm$, then the following are equivalent:
\begin{enumerate}
\item $P$ has TPC.
\item All weight functions for the category of $A_+$-$A_+$bimodules generated by  ${_{A_+}} {\mcal H}_{A_-}$, are trivial.
\item All weight functions for the category of $A_-$-$A_-$ bimodules generated by  ${_{A_+}} {\mcal H}_{A_-}$, are trivial.
\end{enumerate}
\end{lem}
\begin{pf}
Note that by \cite[Remark 6.3]{DGG}, TPC for $P$ is equivalent to all positive weights on $P$ being scalar.
The rest easily follows from Lemma \ref{wtchar}, Proposition \ref{wtext} and Remark \ref{wtextunique}.
\end{pf}

We use the above lemma to prove the following theorem.
\begin{thm}\label{tpc}
If finite index subfactors (of type $II_1$) $N \subset Q$ and $Q \subset M$ have TPC, then $N \subset M$ also has so. More generally, if bifinite bimodules ${_A}{\mcal H}_B$ and ${_B}{\mcal K}_C$ has TPC, then so does ${_A}{\mcal G}_C := {_A}{\mcal H} \us{B}{\otimes} {\mcal K}_C$.
\end{thm}
\begin{pf}
By Lemma \ref{tpceq}, it is enough to prove that all weight functions for the category of $A$-$A$ bimodules generated by ${_A}{\mcal G}_C$, are trivial.
Note that the category of $A$-$A$ bimodules generated by ${_A}{\mcal G}_C$ is same as that by ${_B}{\mcal K} \us{C}{\otimes} \oline{{\mcal G}}_A$.
Again, by lemma \ref{tpceq}, we will be through if we show that any weight function for the category of $B$-$B$ bimodules generated by ${_B}{\mcal K} \us{C}{\otimes} \oline{{\mcal G}}_A$, is trivial.
Now, the $B$-$B$ bimodules generated by ${_B}{\mcal K} \us{C}{\otimes} \oline{{\mcal G}}_A$ are isomorphic to those generated by ${_B}\oline{{\mcal H}}\us{A}{\otimes} {\mcal H}_B$ and ${_B}{\mcal K} \us{C}{\otimes} \oline{{\mcal K}}_B$.
Let $w$ be a weight function for the category of $B$-$B$ bimodules generated by ${_B}{\mcal K} \us{C}{\otimes} \oline{{\mcal K}}_B$
and ${_B}\oline{{\mcal H}}\us{A}{\otimes} {\mcal H}_B$.
Since ${_A}{\mcal H}_B$ (resp., ${_B}{\mcal K}_C$) has TPC, therefore $w$ restricted to the category of $B$-$B$ bimodules generated by ${_B}\oline{{\mcal H}}\us{A}{\otimes} {\mcal H}_B$ (resp.,  ${_B}{\mcal K} \us{C}{\otimes} \oline{{\mcal K}}_B$).
Thus, by the tensor homomorphism property of $w$, we indeed get that $w$ is identically equal to $1$.

The statement involving subfactors can be derived by setting ${_A}{\mcal H}_B = {_N} L^2 (P)_P$ and ${_B}{\mcal K}_C = {_P} L^2 (M)_M$ in the above and using the fact (see proof of \cite[Theorem 5.4]{DGG}) that the normalized bimodule planar algebra associated to ${_N} L^2 (P)_P$ (resp., ${_P} L^2 (M)_M$) is isomorphic to the same associated to $N \subset P$ (resp., $P \subset M$).
\end{pf}
\begin{rem}
In particular, Bisch-Haagerup planar algebras (see \cite{BDG09}) have TPC.
Moreover, given an infinite group $G$ generated by two finite subgroups $H$ and $K$ with nontrivial intersection and acting outerly on a $II_1$-factor $Q$, the planar algebra associated to $Q^H \subset Q \rtimes K$ (which has TPC by  Theorem \ref{tpc} since $Q^H \subset Q$ and $Q \subset Q \rtimes K$ being finite depth / irreducible, have TPC), becomes an example of infinite depth, reducible bimodule planar algebra which can never be perturbed to a non-spherical one.
\end{rem}
\begin{rem}
Given an irreducible, infinite depth, finite index ($>1$) subfactor $N \subset M$, the subfactor $N \subset M_k$ is a reducible and infinite depth; nevertheless, by Theorem \ref{tpc}, its associated planar algebra has TPC, that is, there is no nontrivial weight for all $k \in \N$.
\end{rem}
\section{Perturbation of a bimodule}
If $P$ denotes the bimodule planar algebra coming from a bifinite bimodule ${_A}{\mcal H}_B$ and $z$ is a positive weight on $P$, then by \cite[Theorem 5.13]{DGG}, the perturbation $Q$ of $P$ by $z$, can be associated with a bifinite bimodule ${_C} {\mcal K}_D$ whose bimodule planar algebra is $Q$.
However, there is no direct relation of ${_C} {\mcal K}_D$ with ${_A} {\mcal H}_B$ and the weight $z$.
In this section, we obtain a method of perturbing a bifinite bimodule by a weight to a new one whose associated bimodule planar algebra is the perturbation of the one coming from the bimodule which we start with.
The construction was proposed by Stefaan Vaes; this also answers a question in \cite{DGG}.

We first set up some notations.
Let $R$ denote the hyperfinite $II_{1}$-factor and $\left\{ \theta^{\lambda}\right\} _{\lambda>0}$ be a trace-scaling action of $\mathcal{L}(l^{2}\N)\overline{\otimes}R$, that is, $\theta^{\lambda}\in\mbox{Aut}\left(\mathcal{L}(l^{2}\N)\overline{\otimes}R\right)$ such that $\theta^{\lambda}\circ\theta^{\mu}=\theta^{\lambda\mu}$ and $(\text{Tr}\otimes\tau)\circ\theta^{\lambda}=\lambda(\text{Tr}\otimes\tau)$ where $\tau$ (resp., $\text{Tr}$) is the canonical (resp., semifinite) trace on $R$ (resp., $\mathcal{L}(l^{2}\N)$).
For each $\lambda>0$, set $p^{\lambda}:=\theta^{\lambda}(E_{1,1}\otimes1)\in\mscr P\left(\mathcal{L}(l^{2}\N)\overline{\otimes}R\right)$, $\mathcal{K}^{\lambda}:=p^{\lambda}\left(l^{2}\N\otimes L^{2}(R)\right)$.
Consider the $R$-$R$ bimodule structure on $\mcal K^{\lambda}$ given by
\[
c\cdot p^{\lambda}\left(\xi_{n}\otimes\hat{x}\right)\cdot d:=\theta^{\lambda}(E_{1,1}\otimes c)\left(\xi_{n}\otimes\widehat{xd}\right)
\]
for $x,c,d\in R$, where $\{\xi_{n}\}_{n\in\N}$ denotes the distinguished orthonormal basis of $l^{2} \N$ and $\{E_{i,j}\}_{i,j}$ is a system of matrix units of $\mathcal{L}(l^{2}\N)$ with respect to the basis.
We will now prove a few important properties of these bimodules which will be used in the proof of the main theorem of this section.

First note that $\tau_{\infty}(p^{\lambda})=\lambda$ (where $\tau_{\infty}=\text{Tr}\otimes\tau$) implies $\mbox{dim}\left(\mcal K_{R}^{\lambda}\right)=\lambda$, and since $\mbox{index}\left(\lrsuf R{\mcal K}{\lambda}_{R}\right)=1$, we have $\mbox{dim}\left(\lrsuf R{\mcal K}{\lambda}\right)=\lambda^{-1}$.
Let $\eta_{n}^{\lambda}:=p^{\lambda}\left(\xi_{n}\otimes\hat{1}\right)\in\left(\mcal K^{\lambda}\right)^{o}$ and $\theta_{m,n}^{\lambda}:\mathcal{L}(l^{2}\N)\overline{\otimes}R\ra R$ for $m,n\in\N$ such that $\theta^{\lambda}(\cdot)=\us{m,n}{\sum}E_{m,n}\otimes\theta_{m,n}^{\lambda}(\cdot)$.
Clearly, $p_{m,n}^{\lambda}:=\theta_{m,n}^{\lambda}\left(E_{1,1}\otimes1\right)=\lab\eta_{m}^{\lambda},\eta_{n}^{\lambda}\rab_{R}$.
Thus, $\left\{ \eta_{n}^{\lambda}=\us m{\sum}\xi_{m}\otimes\widehat{p_{m,n}^{\lambda}}\right\} _{n\in\N}$ becomes a basis for $\mcal K_{R}^{\lambda}$.
Also, $\lab\eta_{i}^{\lambda},x\cdot\eta_{j}^{\lambda}\rab_{R}=\theta_{i,j}^{\lambda}(E_{1,1}\otimes x)$ for all $x\in R$, $i,j\in\N$.

\noindent \emph{Standard fact for bimodules:} If $\mcal H$ and $\mcal K$
are $A$-$B$-bimodules with $\{\xi_{i}\}_{i}$ and $\{\zeta_{i}\}_{i}$
as $B$-bases such that $\lab\xi_{i},a\xi_{j}\rab_{B}=\lab\zeta_{i},a\zeta_{j}\rab_{B}$
for all $i,j$ and $a\in A$, then $\xi_{i}\mapsto\zeta_{i}$ extends
to an $A$-$B$ linear unitary.
\begin{lem}\label{udef}
There exists $R$-$R$-linear unitaries $U_{\kappa,\lambda}:\mcal K^{\kappa}\us R{\otimes}\mcal K^{\lambda}\ra\mcal K^{\kappa\lambda}$ and $U_{\kappa,\lambda}^{*}:\mcal K^{\kappa\lambda}\ra\mcal K^{\kappa}\us R{\otimes}\mcal K^{\lambda}$
given by
\[
\eta_{k}^{\kappa}\us R{\otimes}\eta_{l}^{\lambda}\os{U_{\kappa,\lambda}}{\longmapsto}\us i{\sum}\eta_{i}^{\kappa\lambda}\cdot\theta_{i,l}^{\lambda}(E_{k,1}\otimes1)\mbox{ and }\eta_{i}^{\kappa\lambda}\os{U_{\kappa,\lambda}^{*}}{\longmapsto}\us{k,l}{\sum}\eta_{k}^{\kappa}\us R{\otimes}\eta_{l}^{\lambda}\cdot\theta_{l,i}^{\lambda}(E_{1,k}\otimes1).
\]
\end{lem}
\begin{pf}
Since \emph{$\left\{ \eta_{k}^{\kappa}\us R{\otimes}\eta_{l}^{\lambda}:k,l\in\N\right\} $} is a basis for $\mcal K^{\kappa}\us R{\otimes}\mcal K_{R}^{\lambda}$ and $\mbox{dim}\left(\mcal K^{\kappa}\us R{\otimes}\mcal K_{R}^{\lambda}\right)=\kappa\lambda=\mbox{dim}\left(\mcal K_{R}^{\kappa\lambda}\right)$, therefore for existence of $U_{\kappa,\lambda}$, it is enough to prove
\[
\lab\eta_{k}^{\kappa}\us R{\otimes}\eta_{l}^{\lambda},x\cdot\eta_{m}^{\kappa}\us R{\otimes}\eta_{n}^{\lambda}\rab_{R}=\lab U_{\kappa,\lambda}(\eta_{k}^{\kappa}\us R{\otimes}\eta_{l}^{\lambda}),x \cdot U_{\kappa,\lambda}( \eta_{m}^{\kappa}\us R{\otimes}\eta_{n}^{\lambda} ) \rab_{R}.
\]
The right hand side of this equation can be expressed as $\us{i,j}{\sum}\theta_{l,i}^{\lambda}(E_{1,k}\otimes1)\lab\eta_{i}^{\kappa\lambda},x\cdot\eta_{j}^{\kappa\lambda}\rab_{R}\theta_{j,n}^{\lambda}(E_{m,1}\otimes1)$
\begin{align*}
= & \us{i,j}{\sum}\theta_{l,i}^{\lambda}(E_{1,k}\otimes1)\theta_{i,j}^{\kappa\lambda}(E_{1,1}\otimes x)\theta_{j,n}^{\lambda}(E_{m,1}\otimes1)=\us{i,j}{\sum}\theta_{l,i}^{\lambda}(E_{1,k}\otimes1)\theta_{i,j}^{\lambda}\left(\theta^{\kappa}(E_{1,1}\otimes x)\right)\theta_{j,n}^{\lambda}(E_{m,1}\otimes1)\\
= & \theta_{l,n}^{\lambda}\left(E_{1,1}\otimes\theta_{k,m}^{\kappa}(E_{1,1}\otimes x)\right) = \lab \eta_{l}^{\lambda} , \lab \eta_{k}^{\kappa} , x \cdot \eta_{m}^{\kappa} \rab_{R} \cdot \eta_{n}^{\lambda}\rab_{R}=\lab\eta_{k}^{\kappa}\us R{\otimes}\eta_{l}^{\lambda},x\cdot\eta_{m}^{\kappa}\us R{\otimes}\eta_{n}^{\lambda}\rab_{R}.\end{align*}
For $U_{\kappa,\lambda}^{*}$, note that $\lab \eta_{i}^{\kappa\lambda} , x\cdot\eta_{j}^{\kappa\lambda}\rab_{R}=\theta_{i,j}^{\kappa \lambda}(E_{1,1}\otimes x)$ and
\begin{align*}
 & \lab U_{\kappa,\lambda}^{*}(\eta_{i}^{\kappa\lambda}),x\cdot U_{\kappa,\lambda}^{*}(\eta_{j}^{\kappa\lambda})\rab_{R}=\us{k,l,m,n}{\sum}\theta_{i,l}^{\lambda}(E_{k,1}\otimes1)\lab\eta_{k}^{\kappa}\us R{\otimes}\eta_{l}^{\lambda},x\cdot\eta_{m}^{\kappa}\us R{\otimes}\eta_{n}^{\lambda}\rab_{R}\theta_{n,j}^{\lambda}(E_{1,m}\otimes1)\\
= & \us{k,l,m,n}{\sum}\theta_{i,l}^{\lambda}(E_{k,1}\otimes1)\theta_{l,n}^{\lambda}\left(E_{1,1}\otimes\theta_{k,m}^{\kappa}(E_{1,1}\otimes x)\right)\theta_{n,j}^{\lambda}(E_{1,m}\otimes1)=\theta_{i,j}^{\lambda}\left(\theta^{\kappa}(E_{1,1}\otimes x)\right)=\theta_{i,j}^{\kappa \lambda}(E_{1,1}\otimes x).
\end{align*}
Using similar arguments, one can establish the existence of the unitary $U_{\kappa,\lambda}^{*}$. In order to prove that $U_{\kappa,\lambda}^{*}$ (as defined above) is indeed the adjoint of $U_{\kappa,\lambda}$, we consider
\[
U_{\kappa,\lambda}\left(U_{\kappa,\lambda}^{*}(\eta_{i}^{\kappa \lambda})\right)=\us{j,k,l}{\sum} \left[ \eta_{j}^{\kappa\lambda} \cdot \theta_{j,l}^{\lambda}(E_{k,1}\otimes1) \right] \cdot \theta_{li}^{\lambda}(E_{1,k}\otimes1)=\us j{\sum}\eta_{j}^{\kappa\lambda}\cdot\theta_{j,i}^{\lambda}(1)=\eta_{i}^{\kappa\lambda}.
\]
\end{pf}
\begin{lem}
\label{ass}$U_{\kappa,\lambda}$'s are `associative' with respect
to fusion of bimodules, that is, $U_{\kappa\lambda,\mu}\circ(U_{\kappa,\lambda}\us R{\otimes}\mbox{id}_{\mcal K^{\mu}})=U_{\kappa,\lambda\mu}\circ(\mbox{id}_{\mcal K^{\kappa}}\us R{\otimes}U_{\lambda,\mu})$.\end{lem}
\begin{proof}
Note that\begin{eqnarray*}
\eta_{k}^{\kappa}\us R{\otimes}\eta_{l}^{\lambda}\us R{\otimes}\eta_{m}^{\mu} & \os{U_{\kappa,\lambda}\us R{\otimes}\mbox{id}_{\mcal K^{\mu}}}{\longmapsto} & \us i{\sum}\eta_{i}^{\kappa\lambda}\us R{\otimes}\theta_{i,l}^{\lambda}(E_{k,1}\otimes1)\cdot\eta_{m}^{\mu}=\us{i,n}{\sum}\eta_{i}^{\kappa\lambda}\us R{\otimes}\eta_{n}^{\mu}\cdot\theta_{n,m}^{\mu}\left(E_{1,1}\otimes\theta_{i,l}^{\lambda}(E_{k,1}\otimes1)\right)\\
 & \os{U_{\kappa\lambda,\mu}}{\longmapsto} & \us{i,j,n}{\sum}\eta_{j}^{\kappa\lambda\mu}\cdot\theta_{j,n}^{\mu}(E_{i,1}\otimes1)\theta_{n,m}^{\mu}\left(E_{1,1}\otimes\theta_{i,l}^{\lambda}(E_{k,1}\otimes1)\right)\\
 &  & =\us{i,j}{\sum}\eta_{j}^{\kappa\lambda\mu}\cdot\theta_{j,m}^{\mu}\left(E_{i,1}\otimes\theta_{i,l}^{\lambda}(E_{k,1}\otimes1)\right).\end{eqnarray*}
On the other hand,\begin{eqnarray*}
\eta_{k}^{\kappa}\us R{\otimes}\eta_{l}^{\lambda}\us R{\otimes}\eta_{m}^{\mu} & \os{(\mbox{id}_{\mcal K^{\kappa}}\us R{\otimes}U_{\lambda,\mu})}{\longmapsto} & \us s{\sum}\eta_{k}^{\kappa}\us R{\otimes}\eta_{s}^{\lambda\mu}\cdot\theta_{s,m}^{\mu}(E_{l,1}\otimes1)\\
 & \os{U_{\kappa\lambda,\mu}}{\longmapsto} & \us{s,j}{\sum}\eta_{j}^{\kappa\lambda\mu}\cdot\theta_{j,s}^{\lambda\mu}(E_{k,1}\otimes1)\theta_{s,m}^{\mu}(E_{l,1}\otimes1)\end{eqnarray*}
Now, we also have the relation $\us s{\sum}\theta_{j,s}^{\lambda\mu}(E_{k,1}\otimes1)\theta_{s,m}^{\mu}(E_{l,1}\otimes1)=\us s{\sum}\theta_{j,s}^{\mu}\left(\theta^{\lambda}(E_{k,1}\otimes1)\right)\theta_{s,m}^{\mu}(E_{l,1}\otimes1)=\theta_{j,m}^{\mu}\left(\theta^{\lambda}(E_{k,1}\otimes1)(E_{l,1}\otimes1)\right)=\us i{\sum}\theta_{j,m}^{\mu}\left(E_{i,1}\otimes\theta_{i,l}^{\lambda}(E_{k,1}\otimes1)\right).$\end{proof}
\begin{rem}
\label{coherence}By composing the $U_{\lambda,\mu}$'s appropriately
and by Lemma \ref{udef} and \ref{ass}, one can construct a unique
$R$-$R$ linear unitary $U_{\lambda_{1},\ldots,\lambda_{n}}:\mcal K^{\lambda_{1}}\us R{\otimes}\cdots\us R{\otimes}\mcal K^{\lambda_{n}}\ra\mcal K^{\lambda_{1}\cdots\lambda_{n}}$.
\end{rem}
In the next lemma, we will use the following formula:\[
\tau\left(\theta_{i,j}^{\lambda}(E_{k,l}\otimes x)y\right)=\tau_{\infty}\left((E_{j,i}\otimes y)\theta^{\lambda}(E_{k,l}\otimes x)\right)=\lambda\tau_{\infty}\left(\theta^{\lambda^{-1}}(E_{j,i}\otimes y)(E_{k,l}\otimes x)\right)=\lambda\tau\left(\theta_{l,k}^{\lambda^{-1}}(E_{j,i}\otimes y)x\right)\]
for all $x,y\in R$, $i,j,k,l\in\N$.
\begin{lem}
\label{vdef}There exists $R$-$R$-linear unitaries
$V_{\lambda}:\ol{\mcal K^{\lambda}}\ra\mcal K^{\lambda^{-1}}$ and
$V_{\lambda}^{*}:\mcal K^{\lambda^{-1}}\ra\ol{\mcal K^{\lambda}}$
given by\[
\ol{\eta_{i}^{\lambda}}\os{V_{\lambda}}{\mapsto}\lambda^{\frac{1}{2}}\us k{\sum}\eta_{k}^{\lambda^{-1}}\cdot\theta_{k,1}^{\lambda^{-1}}(E_{1,i}\otimes1)\mbox{ and }\eta_{k}^{\lambda^{-1}}\os{V_{\lambda}^{*}}{\mapsto}\lambda^{-\frac{1}{2}}\us i{\sum}\ol{\eta_{i}^{\lambda}\cdot\theta_{i,1}^{\lambda}(E_{1,k}\otimes1)}.\]
\end{lem}
\begin{proof}
We use the above standard fact on unitary between two bimodules. Note
that $\left\{ \ol{\eta_{n}^{\lambda}}\right\} _{n\in\N}$ is a basis
for $\vphantom{\mcal K^{\lambda}}_{R}\ol{\mcal K^{\lambda}}$. For
$x,y\in R$, we have\begin{align*}
\left\langle V_{\lambda}\left(\ol{\eta_{i}^{\lambda}}\right),y\cdot V_{\lambda}\left(\ol{\eta_{j}^{\lambda}}\right)\cdot x\right\rangle  & =\lambda\us{k,l}{\sum}\tau\left(\theta_{1,k}^{\lambda^{-1}}(E_{i,1}\otimes1)\lab\eta_{k}^{\lambda^{-1}},y\cdot\eta_{l}^{\lambda^{-1}}\rab_{R}\theta_{l,1}^{\lambda^{-1}}(E_{1,j}\otimes1)x\right)\\
 & =\lambda\tau\left(\theta_{1,1}^{\lambda^{-1}}(E_{i,j}\otimes y)x\right)=\tau\left(\theta_{j,i}^{\lambda}(E_{1,1}\otimes x)y\right)\end{align*}
which implies $\lsub R{\left\langle V_{\lambda}\left(\ol{\eta_{i}^{\lambda}}\right),V_{\lambda}\left(\ol{\eta_{j}^{\lambda}}\right)\cdot x\right\rangle }=\theta_{j,i}^{\lambda}(E_{1,1}\otimes x)=\lsub R{\left\langle \ol{\eta_{i}^{\lambda}},\ol{\eta_{j}^{\lambda}}\cdot x\right\rangle }$.
Existence of $V_{\lambda}^{*}$ can be obtained using similar arguments
and the basis $\left\{ \eta_{k}^{\lambda^{-1}}\right\} _{k\in\N}$
for $\mcal K_{R}^{\lambda^{-1}}$.
The fact that $V_{\lambda}^{*}$ is indeed the adjoint of  $V_{\lambda}$, comes
from \[
V_{\lambda}\left(V_{\lambda}^{*}\left(\eta_{k}^{\lambda^{-1}}\right)\right)=\us{i,l}{\sum}\theta_{1,i}^{\lambda}(E_{k,1}\otimes1)\cdot\eta_{l}^{\lambda^{-1}}\cdot\theta_{l,1}^{\lambda^{-1}}(E_{1,i}\otimes1)=\us{i,m}{\sum}\xi_{m}\otimes\left[\theta_{m,1}^{\lambda^{-1}}\left(E_{1,i}\otimes\theta_{1,i}^{\lambda}(E_{k,1}\otimes1)\right)\right]^{\widehat{}}\]
and $\us{ i}{\sum}\theta_{m,1}^{\lambda^{-1}}\left(E_{1,i}\otimes\theta_{1,i}^{\lambda}(E_{k,1}\otimes1)\right)=\theta_{m,1}^{\lambda^{-1}}\left(\left(E_{1,1}\otimes1\right)\theta^{\lambda}(E_{k,1}\otimes1)\right)=\us s{\sum}p_{m,s}^{\lambda^{-1}}\theta_{s,1}^{\lambda^{-1}}\left(\theta^{\lambda}(E_{k,1}\otimes1)\right)=p_{m,k}^{\lambda^{-1}}$.\end{proof}
\begin{rem}\label{vinv}
From Lemma \ref{vdef}, it is clear that $\ol{V_{\lambda}(\overline{\xi})}=V_{\lambda^{-1}}^{*}(\xi)$
for all $\xi\in\mcal K^{\lambda}$.
{}
\end{rem}
Next, we consider the $R$-$R$ linear maps \begin{eqnarray*}
 & \mcal K^{\lambda}\us R{\otimes}\ol{\mcal K^{\lambda}}\ni\xi\us R{\otimes}\ol{\eta}\os{F_{\lambda,\mu}}{\longmapsto}\us i{\sum}\lsub R{\lab}\eta,\xi\rab\cdot\alpha_{i}\us R{\otimes}\ol{\alpha_{i}}\in\mcal K^{\mu}\us R{\otimes}\ol{\mcal K^{\mu}}\\
 & \ol{\mcal K^{\lambda}}\us R{\otimes}\mcal K^{\lambda}\ni\ol{\eta}\us R{\otimes}\xi\os{G_{\lambda,\mu}}{\longmapsto}\us j{\sum}\lab\eta,\xi\rab_{R}\cdot\ol{\beta_{j}}\us R{\otimes}\beta_{j}\in\ol{\mcal K^{\mu}}\us R{\otimes}\mcal K^{\mu}\end{eqnarray*}
where $\xi,\eta$ are bounded vectors in $\mcal K^{\lambda}$ and
$\{\alpha_{i}\}_{i}$ and $\{\beta_{j}\}_{j}$ are bases for $\mcal K_{R}^{\mu}$
and $\lrsuf R{\mcal K}{\mu}$ respectively.
\begin{lem}\label{jp}
(i) $F_{\lambda,\mu}=(\lambda\mu)^{\frac{1}{2}}\left(\mbox{id}_{\mcal K^{\mu}}\us R{\otimes}V_{\mu}^{*}\right)\circ U_{\mu,\mu^{-1}}^{*}\circ U_{\lambda,\lambda^{-1}}\circ\left(\mbox{id}_{\mcal K^{\lambda}}\us R{\otimes}V_{\lambda}\right)$,

(ii) $G_{\lambda,\mu}=(\lambda\mu)^{-\frac{1}{2}}\left(V_{\mu}^{*}\us R{\otimes}\mbox{id}_{\mcal K^{\mu}}\right) \circ U_{\mu^{-1},\mu}^{*}\circ U_{\lambda^{-1},\lambda}\circ\left(V_{\lambda}\us R{\otimes}\mbox{id}_{\mcal K^{\lambda}}\right)$.\end{lem}
\begin{proof}
(i) It is equivalent to show that the map\[
\mcal K^{1}\os{U_{\lambda,\lambda^{-1}}^{*}}{\lra}\mcal K^{\lambda}\us R{\otimes}\mcal K^{\lambda^{-1}}\os{\mbox{id}_{\mcal{K^{\lambda}}}\us R{\otimes}V_{\lambda}^{*}}{\lra}\mcal K^{\lambda}\us R{\otimes}\ol{\mcal K^{\lambda}}\os{F_{\lambda,\mu}}{\longmapsto}\mcal K^{\mu}\us R{\otimes}\ol{\mcal K^{\mu}}\os{\mbox{id}_{\mcal K^{\mu}}\us R{\otimes}V_{\mu}}{\longmapsto}\mcal K^{\mu}\us R{\otimes}\mcal K^{\mu^{-1}}\os{U_{\mu,\mu^{-1}}}{\longmapsto}\mcal K^{1}\]
is equal to $(\lambda\mu)^{\frac{1}{2}}\mbox{id}_{\mcal K^{1}}$.
Now,\[
\eta_{1}^{1}\os{U_{\lambda,\lambda^{-1}}^{*}}{\lra}\us{k,l}{\sum}\eta_{k}^{\lambda}\us R{\otimes}\eta_{l}^{\lambda^{-1}}\cdot\theta_{l,1}^{\lambda^{-1}}(E_{1,k}\otimes1)\os{\mbox{id}_{\mcal{K^{\lambda}}}\us R{\otimes}V_{\lambda}^{*}}{\lra}\lambda^{-\frac{1}{2}}\us{i,k,l}{\sum}\eta_{k}^{\lambda}\us R{\otimes}\ol{\eta_{i}^{\lambda}\cdot\theta_{i,1}^{\lambda}(E_{1,l}\otimes1)}\cdot\theta_{l,1}^{\lambda^{-1}}(E_{1,k}\otimes1).\]
We will show that the part $\us{i,k,l}{\sum}\lsub R{\lab}\eta_{i}^{\lambda}\cdot\theta_{i,1}^{\lambda}(E_{1,l}\otimes1),\eta_{k}^{\lambda}\rab\theta_{l,1}^{\lambda^{-1}}(E_{1,k}\otimes1)$
(in the expression obtained by applying $F_{\lambda,\mu}$ on the
above) is equal to $\lambda1$. To see this, note that for all $x\in R$,
we get\[
\us i{\sum}\lab\eta_{i}^{\lambda}\cdot\theta_{i,1}^{\lambda}(E_{1,l}\otimes1),x\cdot\eta_{k}^{\lambda}\rab=\us i{\sum}\tau\left(\theta_{1,i}^{\lambda}(E_{l,1}\otimes1)\lab\eta_{i}^{\lambda},x\cdot\eta_{k}^{\lambda}\rab_{R}\right)=\tau\left(\theta_{1,k}^{\lambda}(E_{l,1}\otimes x)\right)=\lambda\tau\left(\theta_{1,l}^{\lambda^{-1}}(E_{k,1}\otimes1)x\right).\]
So, $\us{i,k,l}{\sum}\lsub R{\lab}\eta_{i}^{\lambda}\cdot\theta_{i,1}^{\lambda}(E_{1,l}\otimes1),\eta_{k}^{\lambda}\rab\theta_{l,1}^{\lambda^{-1}}(E_{1,k}\otimes1)=\lambda\us k{\sum}\theta_{1,1}^{\lambda^{-1}}(E_{k,k}\otimes1)=\lambda1$.
On the other hand, \begin{align*}
\us i{\sum}\eta_{i}^{\mu}\us R{\otimes}\ol{\eta_{i}^{\mu}}\os{\mbox{id}_{\mcal K^{\mu}}\us R{\otimes}V_{\mu}}{\longmapsto} & \mu^{\frac{1}{2}}\us{i,k}{\sum}\eta_{i}^{\mu}\us R{\otimes}\eta_{k}^{\mu^{-1}}\cdot\theta_{k,1}^{\mu^{-1}}(E_{1,i}\otimes1)=\mu^{\frac{1}{2}}U_{\mu,\mu^{-1}}^{*}(\eta_{1}^{1}).\end{align*}

(ii) We follow the same strategy as in (i). So, we consider\[
\eta_{1}^{1}\os{U_{\lambda^{-1},\lambda}^{*}}{\longmapsto}\us{k,l}{\sum}\eta_{k}^{\lambda^{-1}}\us R{\otimes}\eta_{l}^{\lambda}\cdot\theta_{l,1}^{\lambda}(E_{1,k}\otimes1)\os{V_{\lambda}^{*}\us R{\otimes}\mbox{id}_{\mcal{K^{\lambda}}}}{\lra}\lambda^{-\frac{1}{2}}\us{i,k,l}{\sum}\theta_{1,i}^{\lambda}(E_{k,1}\otimes1)\cdot\ol{\eta_{i}^{\lambda}}\us R{\otimes}\eta_{l}^{\lambda}\cdot\theta_{l,1}^{\lambda}(E_{1,k}\otimes1).\]
Now, $\us{i,k,l}{\sum}\theta_{1,i}^{\lambda}(E_{k,1}\otimes1)\cdot\lab\eta_{i}^{\lambda},\eta_{l}^{\lambda}\rab_{R}\cdot\theta_{l,1}^{\lambda}(E_{1,k}\otimes1)=\us k{\sum}\theta_{1,1}^{\lambda}(E_{k,k}\otimes1)=1$.
For other part, note that we need to work with a basis for $\lrsuf R{\mcal K}{\mu}$,
namely $\left\{ V_{\mu^{-1}}\left(\ol{\eta_{i}^{\mu^{-1}}}\right)\right\} _{i\in\N}$
(by Lemma \ref{vdef}). Now, by Remark \ref{vinv} and Lemma \ref{vdef},
we get\[
\us i{\sum}\ol{V_{\mu^{-1}}\left(\ol{\eta_{i}^{\mu^{-1}}}\right)}\us R{\otimes}V_{\mu^{-1}}\left(\ol{\eta_{i}^{\mu^{-1}}}\right)\os{V_{\mu}\us R{\otimes}\mbox{id}_{\mcal K^{\mu}}}{\longmapsto}\mu^{-\frac{1}{2}}\us{i,k}{\sum}\eta_{i}^{\mu^{-1}}\us R{\otimes}\eta_{k}^{\mu}\cdot\theta_{k,1}^{\mu}(E_{1,i}\otimes1)=\mu^{-\frac{1}{2}}U_{\mu^{-1},\mu}^{*}(\eta_{1}^{1}).\]

\end{proof}
\vspace*{1em}

We are now ready to define a `perturbation of a bifinite bimodule'.
Let ${_{A_+}} {\mcal H}_{A_-}$ be a bifinite bimodule for $II_1$-factors $A_\pm$, $P$ be its associated bimodule planar algebra and $z$ be a positive weight on $P$.
The induced weight function for the bicategory of bimodules generated by ${_{A_+}} {\mcal H}_{A_-}$ will be denoted by $w$, and for any $p \in {\mscr P}_{\text{min}} ({\mcal Z} (P_{\eta k}))$ and $p\geq p_0 \in {\mscr P}_{\text{min}} (P_{\eta k})$, we set $w_p := w_{[\text{Range } p_0]}$.
If $A'_{\pm} = R \oline{\otimes} A_{\pm}$, then consider the $A'_{+}$-$A'_{-}$ bimodule
\[
{\mcal H}' := \us{s \in S}{\bigoplus} \; \lsub{ {A'_+} }{ \left[ {\mcal K}^{w_s} \otimes \text{Range } s \right] }_{A'_-}  \text{ where } S = {\mscr P}_{\text{min}} ({\mcal Z} (P_{+1})).
\]
\begin{thm}
The bimodule planar algebra $P'$ associated to $\lsub{A'_+}{{\mcal H}'}_{A'_-}$, is isomorphic to the perturbation $Q$ of $P$ by the weight $z$.
\end{thm}
\begin{pf}
Let $S^k := S \times S \times \dots (k \text{ copies})$ and $\uline{s} \in S^k$ will always denote the $k$ tuple $(s_1, \cdots, s_k)$ for $s_i$'s in $S$.
Consider the bimodules $\lsub{A_+}{ \left( {\mcal H}^{+ \uline{s}} \right) }_{A_{(-)^k}} := \lsub{{A_+}}{ \left[ \text{Range } s_1 \us{A_-}{\otimes} \oline{\text{Range } s_2} \us{A_+}{\otimes} \cdots (k \text{ components}) \right] }_{A_{(-)^k}}$ and  $\lsub{A_-}{ \left( {\mcal H}^{- \uline{s}} \right) }_{A_{(-)^{k-1}}} := \lsub{{A_-}}{ \left[ \oline{\text{Range } s_1} \us{A_+}{\otimes} \text{Range } s_2 \us{A_-}{\otimes} \cdots (k \text{ components}) \right] }_{A_{(-)^{k-1}}}$.
Clearly, ${\mcal H}_{\vlon k}$ is isomorphic (as a bimodule) to $\us{\uline{s} \in S^k}{\bigoplus}  {\mcal H}^{\vlon \uline{s}} $.
So, we may very well assume $P_{\vlon k} = Q_{\vlon k} = {_{A_\vlon}} {\mcal L}_{A_{(-)^k \vlon}} \left( \us{\uline{s} \in S^k}{\bigoplus} {\mcal H}^{\vlon \uline{s}} \right)$ where the action of the tangles are transferred via the bimodule isomorphism.
Similarly, we may also assume
\[
P'_{\vlon k}  = {_{A'_\vlon}} {\mcal L}_{A'_{(-)^k \vlon}} \left( \us{\uline{s} \in S^k}{\bigoplus} \left[ {\mcal K}^{\vlon \uline{s}} \otimes {\mcal H}^{\vlon \uline{s}} \right] \right) \text{ where } {\mcal K}^{\vlon \uline{s}} := \left\{
\begin{array}{ll}
{\mcal K}^{w_{\text{$s_1$}}} \us{R}{\otimes} \oline{ {\mcal K}^{w_{\text{$s_2$}}} } \us{R}{\otimes} \cdots (k \text{ components}) & \text{if } \vlon = +,\\
\oline{ {\mcal K}^{w_{\text{$s_1$}}} } \us{R}{\otimes} {\mcal K}^{w_{\text{$s_2$}}} \us{R}{\otimes} \cdots (k \text{ components}) & \text{if } \vlon = -.
\end{array}
\right.
\]
Set $w_{\vlon \uline{s}} := \ous{k}{\prod}{i=1} w^{(-1)^{i - \delta_{\vlon = +} } }_{s_i}$; from Lemma \ref{udef} and Remark \ref{coherence}, we can obtain an $R$-$R$ linear unitary from ${\mcal K}^{\vlon \uline{s}}$ to ${\mcal K}^{ w_{\vlon \uline{s}} }$ given by
\[
W_{\vlon \uline{s}} = \left\{
\begin{array}{ll}
U_{(w_{s_1},w^{-1}_{s_2},w_{s_3},w^{-1}_{s_4},\ldots)} \circ \left( \text{id}_{\mcal{K}^{w_{s_1}}} \us{R}{\otimes} V_{w_{s_2}} \us{R}{\otimes} \text{id}_{\mcal{K}^{w_{s_3}}} \us{R}{\otimes} V_{w_{s_4}} \us{R}{\otimes} \cdots \right) & \text{if } \vlon = +,\\
U_{(w^{-1}_{s_1},w_{s_2},w^{-1}_{s_3},w_{s_4},\ldots)} \circ \left( V_{w_{s_1}} \us{R}{\otimes} \text{id}_{\mcal{K}^{w_{s_2}}} \us{R}{\otimes} V_{w^{-1}_{s_3}} \us{R}{\otimes} \text{id}_{\mcal{K}^{w_{s_4}}} \cdots \right) & \text{if } \vlon = -.
\end{array}
\right.
\]
We will now try to define an isomorphism between $P'$ and $Q$.
First note that ${_{A_\vlon}} {\mcal L}_{A_{(-)^k \vlon}} \left( {\mcal H}^{\vlon \uline{s}} , {\mcal H}^{\vlon \uline{t}} \right) \neq 0$ if and only if each of ${\mcal H}^{\vlon \uline{s}}$ and ${\mcal H}^{\vlon \uline{t}}$ has at least one irreducible sub-bimodule isomorphic to each other for $\uline{s} , \uline{t} \in S^k$; moreover, under this condition, by the tensor homomorphism property of the weight function $w$, we have $w_{\vlon \uline{s}} = w_{\vlon \uline{t}}$.
Thus, by irreducibility of the $R$-$R$ bimodules ${\mcal K}^\lambda$ for $\lambda > 0$, ${_{A_\vlon}} {\mcal L}_{A_{(-)^k \vlon}} \left( {\mcal H}^{\vlon \uline{s}} , {\mcal H}^{\vlon \uline{t}} \right) \neq 0$ implies and is implied by $\lsub{A'_\vlon}{ {\mcal L} }_{A'_{(-)^k \vlon}} \left( {\mcal K}^{\vlon \uline{s}} \otimes {\mcal H}^{\vlon \uline{s}} \; , \; {\mcal K}^{\vlon \uline{t}} \otimes {\mcal H}^{\vlon \uline{t}} \right) \cong \lsub{R}{ {\mcal L} }_R \left( {\mcal K}^{\vlon \uline{s}}  \; , \; {\mcal K}^{\vlon \uline{t}} \right) \otimes \lsub{A_\vlon}{ {\mcal L} }_{A_{(-)^k \vlon}} \left( {\mcal H}^{\vlon \uline{s}} \; , \; {\mcal H}^{\vlon \uline{t}} \right) \neq 0$;
further the map
\[
{_{A_\vlon}} {\mcal L}_{A_{(-)^k \vlon}} \left( {\mcal H}^{\vlon \uline{t}} , {\mcal H}^{\vlon \uline{s}} \right) \ni x \os{\vphi^{\vlon k}_{\uline{s} , \uline{t}}}{\longmapsto} \left( W^*_{\vlon \ul{s}} \circ W_{\vlon \ul{t}}  \otimes x \right) \in \lsub{A'_\vlon}{ {\mcal L} }_{A'_{(-)^k \vlon}} \left( {\mcal K}^{\vlon \uline{t}} \otimes {\mcal H}^{\vlon \uline{t}} \; , \; {\mcal K}^{\vlon \uline{s}} \otimes {\mcal H}^{\vlon \uline{s}} \right)
\]
is an isomorphism.
This gives rise to a $\ast$-algebra isomorphism $\vphi^{\vlon k} := \left( \left( \vphi^{\vlon k}_{\uline{s} , \uline{t}} \right) \right) : Q_{\vlon k} = P_{\vlon k} \rightarrow P'_{\vlon k}$.

It remains to show that $\vphi$ preserves the action of the tangles $M_{\vlon k}$, $LI_{\vlon k}$, $RI_{\vlon k}$, $E_{\vlon 1}$ for $\vlon \in \{- , +\}$, $k \in \N_0$ (because the action of any other tangle could be obtained from the action of these for a bimodule planar algebra). We are already done with the multiplication tangles. Equivariance under left and right inclusion tangles easily follows from the relations $W^*_{\vlon (\ul{s},s)} \circ W_{\vlon (\ul{t},s)} = \left( W^*_{\vlon \ul{s}} \circ W_{\vlon \ul{t}} \right) \us{R}{\otimes} \text{id}_{{\mcal K}^{w_s} \left/ \ol{{\mcal K}^{w_s}} \right.}$ and $W^*_{-\vlon (s,\ul{s})} \circ W_{-\vlon (s,\ul{t})} = \text{id}_{\left. \ol{{\mcal K}^{w_s}} \right/ {\mcal K}^{w_s} } \us{R}{\otimes} \left( W^*_{\vlon \ul{t}} \circ W_{\vlon \ul{s}} \right)$ whenever $w_{\vlon \ul{s}} = w_{\vlon \ul{t}}$ for $\ul{s} , \ul{t} \in S^k$, $s\in S$.

\noindent {\em Actions of $E_{\pm 1}$:} Note that 
$Q = P^{(z^{1/2} , z^{1/2})}$ implies $C := Q_{E_{+1}} = (z^{1/2} \us{A_-}{\otimes} \text{id}_{\oline{\mcal H}}) \circ P_{E_{+1}} \circ (z^{1/2} \us{A_-}{\otimes} \text{id}_{\oline{\mcal H}})$.
So, by Theorem \ref{bimod-pa-theorem}(i)(d) and the condition $\left. z \right|_{\text{Range } s} = w_s \text{ id}_{\text{Range } s}$ for all $s \in S$, we get
\[
\text{Range } t_1 \us{A_-}{\otimes} \oline{\text{Range } t_2} \ni \xi \us{A_-}{\otimes} \oline{\zeta} \os{ C_{(s_1,s_2) , (t_1,t_2)} }{\longmapsto} \delta_{s_1 = s_2} \sqrt{ w_{s_1} w_{t_1} } \us{i}{\sum} \; \lsub{A_+}{ \langle \zeta , \xi \rangle} \cdot \gamma_i \us{A_-}{\otimes} \oline{\gamma_i} \in \text{Range } s_1 \us{A_-}{\otimes} \oline{\text{Range } s_2} 
\]
where $\xi$ and $\zeta$ are bounded vectors of $\lsub{A_+}{(\text{Range } t_1)}_{A_-}$ and $\lsub{A_+}{(\text{Range } t_2)}_{A_-}$ respectively, and $\{\gamma_i\}_i$ is a basis for $(\text{Range } s_1)_{A_-}$, for $s_1, s_2 , t_1, t_2 \in S$.
Thus, $\left( Q_{E_{+1}} \right)_{(s_1,s_2) , (t_1,t_2)} = \sqrt{ w_{s_1} w_{t_1} } \left( P_{E_{+1}} \right)_{(s_1,s_2) , (t_1,t_2)}$ implying $\left[ \vphi \left( Q_{E_{+1}} \right) \right]_{(s_1,s_2) , (t_1,t_2)} = \sqrt{ w_{s_1} w_{t_1} }  \left( W^*_{+(s_1,s_2)} W_{+(t_1,t_2)} \otimes \left( P_{E_{+1}} \right)_{(s_1,s_2) , (t_1,t_2)} \right)$.
On the other hand, $\left( P'_{E_{+1}} \right)_{(s_1,s_2) , (t_1,t_2)} = F_{t_1,s_1} \otimes  \left( P_{E_{+1}} \right)_{(s_1,s_2) , (t_1,t_2)}$.
Without loss of generality we may assume $s_1 = s_2$ and $t_1=t_2$ because in all other cases,  $\left( P_{E_{+1}} \right)_{(s_1,s_2) , (t_1,t_2)}$, $\left( Q_{E_{+1}} \right)_{(s_1,s_2) , (t_1,t_2)}$, $\left( P'_{E_{+1}} \right)_{(s_1,s_2) , (t_1,t_2)}$ vanish.
Under this assumption and applying Lemma \ref{jp}(i), we get $F_{t_1,s_1} = \sqrt{ w_{s_1} w_{t_1} } W^*_{+(s_1,s_1)} W_{+(t_1,t_1)}$.
Thus, $\vphi \left( Q_{E_{+1}} \right) = P'_{E_{+1}}$.
The same for $E_{-1}$ follows exactly from similar analysis applied on $Q_{E_{-1}} = \left( \text{id}_{\ol{\mcal H}} \us{A_+}{\otimes} z^{-1/2} \right) \circ P_{E_{-1}} \circ \left( \text{id}_{\ol{\mcal H}} \us{A_+}{\otimes} z^{-1/2} \right)$ and Lemma \ref{jp}(ii).

This completes the proof.
\end{pf}
\begin{rem}
An important thing to note is that if the bimodule ${\mcal H}$ is over hyperfinite factors, then the perturbed bimodule ${\mcal H}'$ is also so. This was not the case in the bimodule reconstruction in \cite[Theorem 5.13]{DGG} which used the subfactor reconstruction techniques in \cite{JSW08}.
\end{rem}
\thanks{\noindent
{\bf Acknowledgements.} The authors would like to thank Stefaan Vaes for several valuable discussions, providing us with the main ingredients to answer the questions (1) and (3) mentioned in the introduction, and pointing out  Example \ref{qgpeg}. 
\bibliographystyle{alpha}

\end{document}